\documentclass[a4paper, 11pt]{amsart}

\usepackage{amsmath,amssymb,amsfonts,amsthm}

\usepackage[usenames,dvipsnames, x11names]{xcolor}
\usepackage{csquotes}

\usepackage[all]{xy}
\usepackage[latin1]{inputenc}
\usepackage{enumerate}

\usepackage{xypic}

\usepackage[colorlinks=true,linkcolor=Aquamarine,citecolor=Aquamarine]{hyperref}

\usepackage{tikz}
\usetikzlibrary {patterns,shapes.symbols}

\theoremstyle{theorem}
\newtheorem{theorem}{Theorem}[section]
\newtheorem{conjecture}[theorem]{Conjecture}

\theoremstyle{remark}
\newtheorem{remark}[theorem]{Remark}
\newtheorem{example}[theorem]{Example}

\theoremstyle{definition}
\newtheorem{definition}[theorem]{Definition}

\newcommand{\Mot}{\mathsf{Mot}}
\newcommand{\DM}{\mathsf{DM}}
\newcommand{\DMT}{\mathsf{DMT}}
\newcommand{\MT}{\mathsf{MT}}

\renewcommand{\H}{\mathrm{H}}
\newcommand{\K}{\mathrm{K}}
\newcommand{\M}{\mathrm{M}}
\newcommand{\dR}{\mathrm{dR}}
\newcommand{\B}{\mathrm{B}}
\newcommand{\F}{\mathrm{F}}
\newcommand{\W}{\mathrm{W}}
\newcommand{\gr}{\mathrm{gr}}

\renewcommand{\AA}{\mathbb{A}}
\newcommand{\CC}{\mathbb{C}}
\newcommand{\NN}{\mathbb{N}}
\newcommand{\PP}{\mathbb{P}}
\newcommand{\QQ}{\mathbb{Q}}
\newcommand{\RR}{\mathbb{R}}
\newcommand{\ZZ}{\mathbb{Z}}

\renewcommand{\i}{\mathrm{i}}
\renewcommand{\d}{\mathrm{d}}

\newcommand{\To}{\longrightarrow}

\newcommand{\Li}{\operatorname{Li}}
\renewcommand{\P}{\operatorname{P}}

\newcommand{\comp}{\operatorname{comp}}

\usepackage[all]{xy}
\renewcommand{\diagram}[1]{\SelectTips{cm}{10}\xymatrix{#1}}

\usepackage{tikz}
\usepackage{subfigure}

\title{An introduction to mixed Tate motives}
\author{Cl\'{e}ment Dupont}
\address{Institut Montpelli\'erain Alexander Grothendieck, Universit\'e de Montpellier, CNRS,
Montpellier, France}
\email{clement.dupont@umontpellier.fr}

\begin{document}

\begin{abstract}
Mixed Tate motives are central objects in the study of cohomology groups of algebraic varieties and their arithmetic invariants. They also play a crucial role in a wide variety of questions related to multiple zeta values and polylogarithms, algebraic $\K$-theory, hyperbolic geometry, and particle physics among others. This survey article is an introduction to mixed Tate motives and their many facets. It was written for the proceedings of the \emph{Summer School on Motives and Arithmetic Groups} held in Strasbourg in June 2022.
\end{abstract}

\maketitle

\section{Introduction}

	Motives, as envisioned by Grothendieck, are universal cohomological invariants of algebraic varieties which control some of their arithmetic properties. Grothendieck's first construction of a category of motives was only concerned with the cohomology of smooth projective varieties, and the corresponding motives are now called \emph{pure}. The simplest pure motives are the \emph{pure Tate motives} $\QQ(-n)$, for $n\in\mathbb{Z}$, where $\QQ(-1)$ corresponds to $\H^2(\mathbb{P}^1)$, and $\QQ(-n)=\QQ(-1)^{\otimes n}$, with the usual convention that $\QQ(1)$ is the dual of $\QQ(-1)$.
	
	Motives of general varieties are sometimes called \emph{mixed}, and the notion of weight explains that a mixed motive can be canonically obtained as an iterated extension of pure motives. \emph{Mixed Tate motives} are by definition the iterated extensions of the pure Tate motives. They are very rare among all mixed motives: for instance, the motive corresponding to the $\H^1$ of an algebraic curve is mixed Tate if and only if the curve is rational.
	
	Although pure Tate motives are rather uninteresting in themselves, the study of their iterated extensions reveals an amazing wealth of mathematical structures. The goal of this survey article is to illustrate this fact with several aspects of mixed Tate motives that we now list, referring the reader to the main body of the article for more details and references.\medskip
	
	\paragraph{\textbf{Periods.}} 
	
	Periods of mixed Tate motives include important mathematical constants such as $\pi$, $\log(a)$ for $a\in\QQ_{>0}$, the special values of the Riemann zeta function
	$$\zeta(n)=\sum_{k=1}^\infty \frac{1}{k^n} \quad (n\in\mathbb{N}_{\geq 2}),$$
	and more generally multiple zeta values
	$$\zeta(n_1,\ldots,n_r)=\sum_{1\leq k_1<\cdots <k_r}\frac{1}{k_1^{n_1}\cdots k_r^{n_r}} \quad (n_1,\ldots,n_{r-1}\in \NN_{\geq 1}, n_r\in\NN_{\geq 2}).$$
 	These numbers appear everywhere in particle physics via the computation of Feynman integrals. 
 	
	In the context of mixed Tate motives over a base, periods depend on parameters and we find important special functions such as the classical polylogarithm functions 
	$$\Li_n(z) = \sum_{k=1}^\infty \frac{z^k}{k^n} \quad (n\in\mathbb{N}_{\geq 1}, z\in\CC, |z|<1),$$
	whose role in the computation of volumes in hyperbolic geometry was discovered by Lobachevsky in the 1830s.\medskip
	
	\paragraph{\textbf{Algebraic $\K$-theory.}}
	
	It is the study of invariants $\K_i(R)$, for $i\in\ZZ$, associated to a ring $R$. They are abelian groups that should be thought of as homotopical invariants of the category of finitely generated projective $R$-modules (also known as vector bundles on $\operatorname{Spec}R$). The abelian category $\MT(F)$ of mixed Tate motives over a field $F$ is expected to be related to the rational $\K$-theory of $F$ via the isomorphism
	\begin{equation}\label{eq: Beilinson formula}
	\operatorname{Ext}^i_{\MT(F)}(\QQ(-n),\QQ(0)) \stackrel{?}{\simeq} \gr_\gamma^n \K_{2n-i}(F)_\QQ.
	\end{equation}
	Here $\gr_\gamma^n$ denotes the $n$-th graded piece of the $\gamma$-filtration in $\K$-theory, and $V_\QQ$ denotes the rationalization $V\otimes_\ZZ\QQ$ of an abelian group $V$. The groups \eqref{eq: Beilinson formula} are referred to in the literature as (rational) \emph{motivic cohomology groups} of the field $F$. The existence of the abelian category of mixed Tate motives over a field, and the isomorphism \eqref{eq: Beilinson formula}, are conjectural in general but both are known in the case of number fields. \medskip
	
	\paragraph{\textbf{Dedekind zeta values.}}
	
	Recall the Dedekind zeta function of a number field $F$,
	$$\zeta_F(s) = \sum_{\mathfrak{a}}\operatorname{N}(\mathfrak{a})^{-s} \quad (s\in\CC,\operatorname{Re}(s)>1),$$
	where the sum ranges over the non-zero ideals of the ring of integers $\mathcal{O}_F$, and $\operatorname{N}(\mathfrak{a})=|\mathcal{O}_F/\mathfrak{a}|$ is the norm. The analytic class number formula expresses its residue at $s=1$ in terms of important arithmetic invariants of $F$, including a transcendental quantity, the \emph{regulator}, which is a determinant of logarithms of units of $\mathcal{O}_F$. More generally, a theorem of Borel relates the special value $\zeta_F(n)$, for an integer $n\geq 2$, to a \emph{higher regulator} defined on $\K_{2n-1}(F)$.
	
	The explicit computation of higher regulators is difficult and requires a fine understanding of the structure of the category of mixed Tate motives, via \eqref{eq: Beilinson formula}. For instance, Zagier's conjecture predicts that $\zeta_F(n)$ can be expressed in terms of special values of the $n$-th polylogarithm function $\Li_n$ at elements of $F$.\medskip
	
	\paragraph{\textbf{Notation and conventions.}}	
	
	Throughout this article we use the notation $V_\QQ:=V\otimes_\ZZ\QQ$ for $V$ an abelian group and $V_\CC:=V\otimes_\QQ\CC$ for $V$ a $\QQ$-vector space.\medskip
	
	\paragraph{\textbf{Acknowledgements.}} This survey paper is based on a mini-course given in Strasbourg in June 2022 on the occasion of the IRMA Summer School on Motives and Arithmetic Groups, and I would like to thank the organizers and all the participants of the school for their interesting questions and feedback. Many thanks to Fran\c{c}ois Fillastre for stimulating discussions on hyperbolic geometry, and to Javier Fres\'{a}n for his many comments on a first version of this paper.

\section{Cohomology and periods}

	We review some (classical and not-so-classical) facts and constructions on Betti and de Rham cohomology of algebraic varieties, and periods. In passing we study families of examples which we will later lift to mixed Tate motives.

	\subsection{Betti and de Rham}\label{par: betti de rham}
		
			We will be interested in two incarnations of the cohomology of an algebraic variety $X$ defined over a field $F$.
			\begin{enumerate}[$\bullet$]
			\item \emph{Algebraic de Rham cohomology} is available if $F$ has characteristic zero and produces finite dimensional vector spaces over $F$ denoted by $\H^n_\dR(X)$. For a smooth variety $X$ they were defined by Grothendieck \cite{grothendieckderham} as the hypercohomology groups (in the Zariski topology) of the complex of algebraic differential forms:
			$$\H^n_\dR(X):=\mathbb{H}^n(X,\Omega^\bullet_{X/F}).$$
			We refer the reader to \cite[Chapter 3]{hubermuellerstach} for the case of singular varieties.
			\item \emph{Betti cohomology} is available if $F$ has an embedding $\sigma\colon F\to \CC$ and produces finite dimensional vector spaces over $\QQ$ denoted by $\H^n_{\B,\sigma}(X)$, or $\H^n_\B(X)$ if $\sigma$ is understood. They are defined as the singular cohomology groups  with rational coefficients of the analytic variety~$X_\sigma^{\mathrm{an}} := (X\times_{F,\sigma}\CC)^{\mathrm{an}}$. In other words, they are dual to the singular homology groups of $X_\sigma^{\mathrm{an}}$ with rational coefficients:
			$$\H^n_{\B,\sigma}(X) := \H^n_{\mathrm{sing}}(X_\sigma^{\mathrm{an}};\QQ)=\H_n^{\mathrm{sing}}(X_\sigma^{\mathrm{an}};\QQ)^\vee.$$
			\end{enumerate}
			By de Rham \cite{derhamthese} and Grothendieck \cite{grothendieckderham}, integration of differential forms on cycles induces a canonical $\CC$-linear comparison isomorphism
			\begin{equation}\label{eq: comparison isomorphism}
			\operatorname{comp}_\sigma\colon \H^n_\dR(X)\otimes_{\,F,\sigma}\CC \stackrel{\sim}{\To} \H^n_{\B,\sigma}(X)\otimes_{\, \QQ}\CC.
			\end{equation}
			Its matrix, relative to the choice of an $F$-basis of $\H^n_\dR(X)$ and a $\QQ$-basis of $\H^n_{\B,\sigma}(X)$, is called a \emph{period matrix} of $\H^n(X)$. (We will use the symbol $\H^n(X)$ for the package consisting of de Rham cohomology and Betti cohomology of $X$ together with the comparison isomorphism.)

			\begin{example}\label{ex: pi}
			For $F=\QQ$, a period matrix of the cohomology group $\H^1(\AA^1\setminus \{0\})$ is the $1\times 1$ matrix
			\begin{equation}\label{eq: period matrix pi}
			\left(\begin{matrix} \,2\pi \i\,\end{matrix} \right).
			\end{equation}
			Indeed, a basis of $\H^1_\dR(\AA^1\setminus\{0\})$ is given by the class of the differential form $\frac{\d x}{x}$, a basis of $\H^1_\B(\AA^1\setminus \{0\})^\vee = \H_1^{\mathrm{sing}}(\CC^*;\QQ)$ is given by the class of the loop $\gamma\colon t\mapsto e^{2\pi \i t}$ around $0$, and the corresponding period is
			$$\int_\gamma\frac{\d x}{x} = 2\pi\i.$$
			\end{example}
	
			We will also consider the \emph{relative cohomology group} $\H^n(X,Y)$ where $Y$ is a closed subvariety of $X$. In the de Rham setting, if both $X$ and $Y$ are smooth, it is defined as the hypercohomology of the shifted cone of the restriction map~$\Omega^\bullet_{X/F}\to i_*\Omega^\bullet_{Y/F}$, where $i\colon Y\hookrightarrow X$ denotes the closed immersion. More precisely,
			$$\Omega^\bullet_{(X,Y)/F} = \Omega^\bullet_{X/F}\oplus i_*\Omega^{\bullet-1}_{Y/F}$$
			where the differential is the sum of the de Rham differentials on $X$ and $Y$ and the restriction, with appropriate signs. We refer the reader to \cite[Chapter 3]{hubermuellerstach} for the general case. In the Betti setting, it is defined as the relative singular cohomology of the pair $(X_\sigma^{\mathrm{an}},Y_\sigma^{\mathrm{an}})$.
			
			\begin{example}\label{ex: log}
			Let $F$ be a subfield of $\CC$ and fix $z\in F\setminus\{0,1\}$. The relative cohomology group $\H^1(\mathbb{A}^1\setminus \{0\},\{1,z\})$ has dimension $2$, as can be seen from the long exact sequence in relative cohomology:
						$$0\to \H^0(\mathbb{A}^1\setminus \{0\})\to \H^0(\{1,z\}) \to \H^1(\mathbb{A}^1\setminus \{0\},\{1,z\}) \to \H^1(\mathbb{A}^1\setminus \{0\}) \to 0.$$
			With well-chosen bases we get the period matrix
			\begin{equation}\label{eq: period matrix log}
			\setlength{\arraycolsep}{3pt}\def\arraystretch{1.3}
			\left( \begin{matrix} 1 & \log(z) \\ 0 & 2\pi \i \end{matrix} \right),
			\end{equation}
			which features (a determination of) the logarithm of $z$, defined as an integral over a path from $1$ to $z$ in $\CC^*$:
			$$\int_1^z\frac{\d x}{x} = \log(z).$$
			\end{example}
			
			If $F$ is algebraic over $\QQ$, the entries of period matrices of relative cohomology groups span a subalgebra of $\CC$ which coincides with the algebra of effective periods in the sense of Kontsevich--Zagier \cite{kontsevichzagier}, as proved in \cite[Theorem 12.2.1]{hubermuellerstach}.
			
	\subsection{Cohomology of families of algebraic varieties}\label{par: families}
	
		Algebraic varieties often come in families, i.e., as morphisms $\pi\colon X\to S$. We first focus on the case where $X$ and $S$ are smooth over a subfield $F$ of $\CC$, and $\pi$ is smooth and proper. This ensures, via Ehresmann's theorem, that $\pi^{\mathrm{an}}:X^{\mathrm{an}}\to S^{\mathrm{an}}$ is a locally trivial fibration.
		\begin{enumerate}[$\bullet$]
		\item The algebraic de Rham cohomology groups of the fibers of $\pi$ assemble into an algebraic vector bundle on $S$:
		$$\mathcal{H}^n_\dR(X/S) := \mathbb{R}^n\pi_*(\Omega^\bullet_{X/S}).$$
		It comes equipped with a flat algebraic connection
		$$\nabla\colon \mathcal{H}^n_\dR(X/S) \To   \mathcal{H}^n_\dR(X/S)\otimes_{\mathcal{O}_S} \Omega^1_{S/F} $$
		called the \emph{Gauss--Manin} connection.
		\item The Betti cohomology groups of the fibers of $\pi$ assemble into a local system of finite dimensional vector spaces over $\QQ$ on $S^{\mathrm{an}}$:
		$$\mathcal{H}^n_{\B}(X/S):=\operatorname{R}^n\!\pi^{\mathrm{an}}_*(\QQ_{X^{\mathrm{an}}}).$$
		If $S^{\mathrm{an}}$ is connected and equipped with a basepoint $s_0$, this local system is completely determined by the \emph{monodromy representation}
		$$\pi_1(S^{\mathrm{an}},s_0) \To \operatorname{Aut}_\QQ(\H^n_\B(X_{s_0})).$$
			\end{enumerate}
			
			The comparison isomorphisms \eqref{eq: comparison isomorphism} for the fibers of $\pi$ assemble into an isomorphism of analytic vector bundles with flat connection on $S^{\mathrm{an}}$,
			$$\left((\mathcal{H}^n_\dR(X/S)\otimes_{\mathcal{O}_S}\mathcal{O}_{S_\CC})^{\mathrm{an}},\nabla^{\mathrm{an}}\right) \stackrel{\sim}{\To} \left(\mathcal{H}^n_\B(X/S) \otimes_{\,\QQ_{S^{\mathrm{an}}}}\mathcal{O}_{S^{\mathrm{an}}},\operatorname{id}\otimes\, \d\right),$$
			where $\d:\mathcal{O}_{S^{\mathrm{an}}}\to \Omega^1_{S^{\mathrm{an}}}$ is the exterior differential on holomorphic forms. Equivalently we get an isomorphism of complex local systems
			$$\left((\mathcal{H}^n_\dR(X/S)\otimes_{\mathcal{O}_S}\mathcal{O}_{S_\CC})^{\mathrm{an}}\right)^{\nabla^{\mathrm{an}}} \stackrel{\sim}{\To} \mathcal{H}^n_\B(X/S)\otimes_{\,\QQ_{S^{\mathrm{an}}}} \CC_{S^{\mathrm{an}}}.$$
			The entries of period matrices are multi-valued holomorphic functions on $S^{\mathrm{an}}$ and are sometimes called \emph{period functions}. They naturally form solutions to a first-order algebraic system of linear differential equations on $S$, which is simply the equation $\nabla=0$ in some local trivialization of $\mathcal{H}^n_\dR(X/S)$.
			
			This discussion generalizes to the relative cohomology of a pair $(X,Y)$ over $S$ if the pair $(X^{\mathrm{an}},Y^{\mathrm{an}})\to S^{\mathrm{an}}$ is locally trivial, i.e., locally a product of $S^{\mathrm{an}}$ with a pair of varieties. The generalization to more general situations, and in particular to non-smooth morphisms $\pi$, requires the formalism of algebraic $\mathcal{D}$-modules (on $S$) and constructible sheaves (on $S^{\mathrm{an}}$).
			
			\begin{example}\label{ex: log variation}
			One may recast Example \ref{ex: log} in the language of (relative) cohomology of algebraic varieties over the base $S=\mathbb{A}^1_\QQ\setminus \{0\}$, which is the space of parameters $z$. (Strictly speaking, Example \ref{ex: log} is attached to the family $(\mathbb{A}^1_\QQ\setminus\{0\}, \{1,z\})$ parametrized by a point $z\in\mathbb{A}^1_\QQ\setminus\{0,1\}$, which naturally extends to $S$.)
			In the de Rham setting we obtain the trivial vector bundle $\mathcal{O}_S^{2}$ with connection matrix
			$$\setlength{\arraycolsep}{4pt}\def\arraystretch{1.2}\omega =\left(\begin{matrix}0 & \frac{\d z}{z} \\ 0 & 0\end{matrix}\right).$$
			This corresponds to the fact that the rows of the period matrix \eqref{eq: period matrix log} satisfy the differential equation $\d f=f\omega$. In the Betti setting we obtain a rank $2$ local system on $S^{\mathrm{an}}=\mathbb{C}^*$, which has fiber $\QQ^2$ at the basepoint $1$ and monodromy matrix
			$$\setlength{\arraycolsep}{4pt}\mu = \left( \begin{matrix} 1 & 1 \\ 0 & 1 \end{matrix}\right).$$
			This reflects the fact that the columns of the period matrix \eqref{eq: period matrix log} are multi-valued functions on $\CC^*$ which transform as $g\leadsto \mu g$ when $z$ winds positively around $0$.
			\end{example}
		
		\subsection{Classical polylogarithms as period functions}
		
		Classical polylogarithms are a natural generalization of the logarithm function (see the nice survey \cite{hainclassical}). For an integer $n\geq 1$, the $n$-th polylogarithm function is defined by the power series
		$$\Li_n(z) = \sum_{k=1}^\infty \frac{z^k}{k^n} \quad (z\in\CC, |z|<1).$$
		We have $\Li_1(z) = -\log(1-z)$, and the recursive relation
		\begin{equation}\label{eq: diff eq Li n}
		\d\Li_n(z) = \Li_{n-1}(z) \frac{\d z}{z}
		\end{equation}
		for $n\geq 2$, which explains that every $\Li_n$ extends to a multi-valued holomorphic function on $\CC\setminus \{1\}$. Let us introduce the $(n+1)\times (n+1)$ matrix
		\begin{equation}\label{eq: period matrix polylog}
		\setlength{\arraycolsep}{5pt} \def\arraystretch{1.4}
			 \left(\begin{matrix} 
			1 & \Li_1(z) & \Li_2(z)        & \Li_3(z) & \cdots & \Li_n(z)\\
			 & 2\pi \i   & 2\pi\i\log(z)  & 2\pi \i \frac{\log^2(z)}{2}& \cdots & 2\pi \i\frac{\log^{n-1}(z)}{(n-1)!}\\
			 &            & (2\pi \i)^2      & (2\pi \i)^2 \log(z) & \cdots & (2\pi \i)^2 \frac{\log^{n-2}(z)}{(n-2)!} \\
			 &  &  & (2\pi \i)^3 & &  \\
			 & & & & & \vdots \\
			& & 0 & & \ddots &  \\
			& & & & & \\
			& & & & &  (2\pi \i)^n 
			\end{matrix}\right)
			\end{equation}
			whose entries are viewed as multi-valued holomorphic functions on $\CC\setminus \{0,1\}$. The rows of \eqref{eq: period matrix polylog} are a system of fundamental solutions of the differential equation $\d f=f\omega$ where 
			$$\setlength{\arraycolsep}{5pt}\def\arraystretch{1.3}
			\omega = \left(\begin{matrix}
			0 & \frac{\d z}{1-z} &      &  &&\\
			 & 0   & \frac{\d z}{z}&  &  0 & \\
			 &      & 0     & \frac{\d z}{z} && \\
			 &  &  & 0  &\ddots &  \\
			 &0 & & &\ddots &\frac{\d z}{z} \\
			& &  & & &  0 \\
			\end{matrix}\right).$$
			When $z$ winds positively around $0$ and $1$, the columns of \eqref{eq: period matrix polylog} transform as $g\leadsto \mu_0 g$ and $g\leadsto \mu_1 g$ respectively, where
		$$
		\setlength{\arraycolsep}{5pt}\def\arraystretch{1.3}
		\mu_0 = \left(\begin{matrix}
		1 & 0 & 0  & 0 & & 0\\
		& 1 & 1  &\frac{1}{2}&& \frac{1}{(n-1)!} \\
		&&1 &1&&\\
		&0&&\ddots&\ddots&\frac{1}{2}\\
		&&&&  1& 1\\
		&&&&& 1
		\end{matrix}\right)
		\quad \mbox{ and } \quad 
		\mu_1 = \left(\begin{matrix}  
		1 & -1 &  && \\
		0 & \,1 &  &0& \\
		&&1 &&\\
		&0&&\ddots&\\
		&&&&  1
		\end{matrix}\right).$$

\begin{figure}
		\begin{center}
		\begin{tikzpicture}
		\def\a{2.5};
		\draw[thick,VioletRed3] (0,-.3*\a) -- (0,2.2*\a);
		\draw[thick,VioletRed3] (\a,-.3*\a) -- (\a,2.2*\a);
		\draw[thick,VioletRed3] (-.3*\a,0) -- (2.2*\a,0);
		\draw[thick,VioletRed3] (-.3*\a,\a) -- (2.2*\a,\a);
		\draw[thick,VioletRed3,fill=VioletRed3!20!white] (0,0) -- (0,\a) -- (\a,\a) -- (\a,0) -- cycle;
		\path[draw, color=Turquoise4, thick] (.7*\a, 2.2*\a) to[out = -80, in = 170] (2.2*\a,.7*\a);
		\end{tikzpicture}
		\end{center}	
		\caption{The hypersurfaces $\{zx_1x_2=1\}$ and $\{x_1(1-x_1)x_2(1-x_2)=0\}$ in the affine plane.}
		\label{fig: polylog}
		\end{figure}
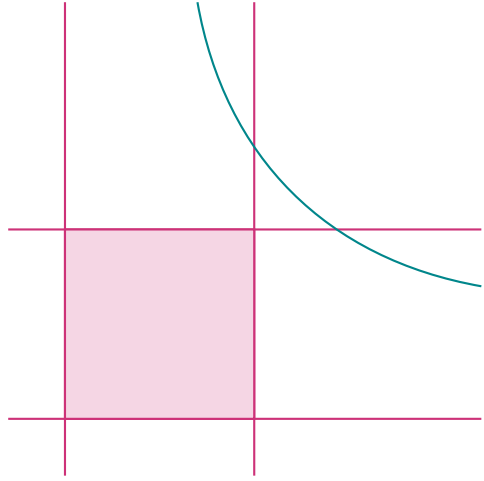

		It is explained in \cite{dupontfresan} that \eqref{eq: period matrix polylog} is the period matrix of a family of relative cohomology groups on $S=\mathbb{A}^1_\QQ\setminus \{0,1\}$ whose fiber at $z$ is 
		\begin{equation}\label{eq: relative cohomology polylog}
		\H^n(\mathbb{A}^n\setminus \{z x_1\cdots x_n=1\}, \{x_1(1-x_1)\cdots x_n(1-x_n)=0\}).
		\end{equation}
		These cohomology groups naturally stem from the integral formula
		\begin{equation}\label{eq: polylog integral}
		\Li_n(z) = \int_{[0,1]^n} \frac{z\,\d x_1\cdots \d x_n}{1-z x_1\cdots x_n},
		\end{equation}
		valid for instance if $z\notin (1,\infty)$. Indeed, the differential form in \eqref{eq: polylog integral} has poles along the hypersurface $\{z x_1\cdots x_n=1\}$, while the integration domain has its boundary contained in the complex points of the hypersurface ${\{x_1(1-x_1)\cdots x_n(1-x_n)=0\}}$. Figure \ref{fig: polylog} illustrates the corresponding geometry in the case $n=2$.
		
	\subsection{Single-valued periods}
		
		Let us consider an embedding $\sigma\colon F\to \CC$ and its complex conjugate $\overline{\sigma}\colon F\to \CC$. For an algebraic variety $X$ defined over $F$, complex conjugation gives rise to an anti-holomorphic isomorphism $X_{\sigma}^{\mathrm{an}} \To X_{\overline{\sigma}}^{\mathrm{an}}$, which induces a $\QQ$-linear isomorphism
		$$\operatorname{F}_\infty\colon \H^n_{\B,\sigma}(X) \stackrel{\sim}{\To} \H^n_{\B,\overline{\sigma}}(X)$$
		called the \emph{Frobenius at infinity} associated to the pair $(\sigma,\overline{\sigma})$. Following Brown \cite{brownsvmzv} we define the isomorphism
		\begin{equation}\label{eq: sv sigma general}
		\operatorname{s}_\sigma\colon \H^n_\dR(X)\otimes_{F,\sigma}\CC \stackrel{\sim}{\To} \H^n_\dR(X)\otimes_{F,\overline{\sigma}}\CC
		\end{equation}
		by the following commutative diagram:
		$$
		\xymatrixcolsep{1cm}\xymatrixrowsep{1cm}\begin{diagram}
		{
		\H^n_\dR(X)\otimes_{F,\sigma}\CC  \ar[r]^-{\operatorname{comp}_\sigma} \ar[d]_{\operatorname{s}_\sigma} & \H^n_{\B,\sigma}(X)\otimes_\QQ\CC  \ar[d]^-{\operatorname{F}_\infty\otimes\operatorname{id}} \\ 
		\H^n_\dR(X)\otimes_{F,\overline{\sigma}}\CC & \H^n_{\B,\overline{\sigma}}(X)\otimes_\QQ\CC   \ar[l]^-{\operatorname{comp}^{-1}_{\overline{\sigma}}}   &
		}
		\end{diagram}
		$$
		It is called the \emph{single-valued period isomorphism} associated to the pair $(\sigma,\overline{\sigma})$. Its coefficients are called \emph{single-valued periods} for $\H^n(X)$.
		
		The situation is simpler if $\sigma$ is a real embedding, in which case we view $F\subset \RR$ and drop $\sigma$ from the notation. The single-valued period isomorphism is then defined over $\RR$, i.e., induces an $\RR$-linear involution
		$$\operatorname{s}\colon \H^n_\dR(X)\otimes_F \RR \stackrel{\sim}{\To} \H^n_\dR(X)\otimes_F\RR.$$
		If $P$ denotes a period matrix for $\H^n(X)$ then a matrix for $\operatorname{s}$ is given by the product
		$$\overline{P}^{\, -1}P,$$
		which we call a \emph{single-valued period matrix} for $\H^n(X)$.
		
		\begin{example}\label{ex: sv pi}
		In the setting of Example \ref{ex: pi}, a period matrix $P$ is \eqref{eq: period matrix pi} and a single-valued period matrix is the $1\times 1$ matrix
		\begin{equation}
		\overline{P}^{\, -1}P=\left(\begin{matrix} \,-1\,\end{matrix} \right).
		\end{equation}
		Therefore the ``single-valued version'' of $2\pi \i$ is $-1$.
		\end{example}
		
		The term ``single-valued period'' is explained by the general setting of families of algebraic varieties (see \S\ref{par: families}): the entries of period matrices are multi-valued holomorphic functions on $S^{\mathrm{an}}$, whereas the entries of single-valued period matrices are \emph{single-valued} real-analytic functions. This is because the single-valued period isomorphisms assemble into an endomorphism of the analytic vector bundle $(\mathcal{H}^n_\dR\otimes_{\mathcal{O}_S}\mathcal{O}_{S_\CC})^{\mathrm{an}}$. This can also be viewed more concretely in matrix form: by following a loop in $S^{\mathrm{an}}$ the period matrix changes as $P\leadsto \mu P$ where $\mu$ is a monodromy matrix, and the entries of the single-valued period matrix do not change:
		$$\overline{P}^{\, -1} P \leadsto \overline{\mu P}^{\, -1}\mu P = \overline{P}^{\,-1}\overline{\mu}^{\, -1}\mu P = \overline{P}^{\, -1}P$$
		because $\mu$ has rational entries. 
		
		\begin{example}\label{ex: sv log}
		In the setting of Example \ref{ex: log}, taking for $\sigma$ the implicit embedding of $F$ inside $\CC$, a period matrix $P$ is \eqref{eq: period matrix log} and a single-valued period matrix is
		$$\setlength{\arraycolsep}{3pt}\def\arraystretch{1.3}
		\overline{P}^{\, -1}P = \left(\begin{matrix} 1 & 2\log|z| \\ 0 & -1\end{matrix}\right),$$
		which features the real-analytic ``single-valued version'' $2\log|z|$ of the multi-valued holomorphic function $\log(z)$.
		\end{example}
		
		\begin{remark}
		Computing single-valued periods is an \emph{a priori} difficult task since one needs to know an entire period matrix in order to compute any given entry of a single-valued period matrix. However, one can give explicit integral formulas for single-valued periods in certain natural geometric situations \cite{BD1} (see also \cite{BD2} for an application to periods arising in string theory). For instance, in the setting of Example \ref{ex: sv log} we get the following integral formula for the ``single-valued logarithm'':
		$$2\log|z| = \frac{1}{2\pi\i}\iint_{\PP^1(\CC)} \operatorname{dlog}\left(\frac{x-z}{x-1}\right)\wedge \frac{\d\overline{x}}{\overline{x}}.$$
		\end{remark}
		
\subsection{Single-valued polylogarithms as single-valued period functions}\label{subsec: sv polylog}

		Let us apply the formalism of single-valued periods to the case of classical polylogarithms, following Beilinson--Deligne \cite{beilinsondeligne}. We first focus on the case of the dilogarithm and the period matrix 
		$$\setlength{\arraycolsep}{3pt}\def\arraystretch{1.3}
		P = \left(\begin{matrix} 1 & - \log(1-z) & \Li_2(z) \\0 & 2\pi\i & 2\pi\i\log(z) \\ 0 & 0 & (2\pi \i)^2\end{matrix}\right),$$
		which is the case $n=2$ of \eqref{eq: period matrix polylog}. The corresponding single-valued period matrix is
		$$\setlength{\arraycolsep}{3pt}\def\arraystretch{1.3}
		\overline{P}^{\, -1}P = \left(\begin{matrix} 1 & -2\log|1-z| & 2\i\operatorname{Im}(\Li_2(z))-2\log|z|\overline{\log(1-z)} \\0 & -1 & -2\log|z| \\ 0 & 0 & 1\end{matrix}\right).$$
		For reasons that will become clear later (\S\ref{subsec: regulators sv}), it is natural to consider instead the half logarithm of the unipotent matrix $D\overline{P}^{\, -1}P$, where $D$ denotes the diagonal matrix with entries $1,-1,1$:
		$$\setlength{\arraycolsep}{4pt}\def\arraystretch{1.3}
		\frac{1}{2}\log(D\overline{P}^{\, -1}P) = \left(\begin{matrix} 0 & -\log|1-z| & \i\operatorname{Im}(\Li_2(z)+\log|z|\log(1-z)) \\0 & 0 & \log|z| \\ 0 & 0 & 0\end{matrix}\right).$$
		We encounter an important special function called the \emph{Bloch--Wigner dilogarithm},
		$$\P_2(z) := \operatorname{Im}(\Li_2(z)+\log|z|\log(1-z)).$$
		It is a single-valued real-analytic function on $\CC\setminus\{0,1\}$, which extends continuously to $\mathbb{P}^1(\CC)$ by setting $\P_2(0)=\P_2(1)=\P_2(\infty)=0$.
		
		By performing the same operation on the period matrix \eqref{eq: period matrix polylog} for the classical polylogarithm $\Li_n(z)$, we discover a well-behaved ``single-valued polylogarithm''
		\begin{align*}
		\operatorname{P}_n(z) & := \operatorname{Re}_n\left(\sum_{k=0}^{n-1} \frac{2^kB_k}{k!}\log^k|z|\Li_{n-k}(z)\right) \\
		& = \operatorname{Re}_n(\Li_n(z) - \log|z|\Li_{n-1}(z)+\cdots),
		\end{align*}
		where $\operatorname{Re}_n$ is the real part $\operatorname{Re}$ if $n$ is odd and the imaginary part $\operatorname{Im}$ if $n$ is even, and $B_k$ are the Bernoulli numbers defined by their exponential generating series
		$$\sum_{k=0}^\infty B_k\frac{x^k}{k!} = \frac{x}{e^x-1}.$$

\section{First intermezzo: mixed Hodge theory}
	
		The Betti cohomology groups of a complex algebraic variety are canonically endowed with a linear algebra datum called a mixed Hodge structure, which we now discuss. Our goal here is to introduce mixed Hodge--Tate structures, the iterated extensions of the pure Hodge--Tate structures $\QQ(-n)$, which are realizations of mixed Tate motives and share several formal features with them.
		
	\subsection{Pure Hodge structures}
	
		\begin{definition}
		Let $w\in\ZZ$. A \emph{pure Hodge structure} of weight $w$ is the datum of a finite dimensional vector space $H$ over $\QQ$ and a $\CC$-linear decomposition
		$$H_\CC = \bigoplus_{p+q=w} H^{p,q}$$
		where the indices $p,q$ are integers, called the \emph{Hodge decomposition}, which satisfies the \emph{Hodge symmetry}
		$$\overline{H^{p,q}}=H^{q,p} \quad \mbox{ for all } p,q.$$
		Here $\overline{(\cdot)}$ denotes the complex conjugation operator on $H_\CC$. A \emph{morphism of pure Hodge structures of weight $w$} is a $\QQ$-linear map whose complexification is compatible with the Hodge decompositions.
		\end{definition}
	
		Note that a pure Hodge structure of odd weight must have even dimension. This explains why the simplest pure Hodge structures have even weight.
		
		\begin{definition}
		The \emph{pure Hodge--Tate structure} $\QQ(-n)$ has underlying vector space $H=\QQ$ and Hodge decomposition $H_\CC=\CC=H^{n,n}$. It is up to isomorphism the unique pure Hodge structure of weight $2n$ and dimension $1$.
		\end{definition}
		
		The importance of Hodge structures in the study of the topology of algebraic varieties comes from the following landmark result of Hodge \cite{hodge}, for which the reader can refer to \cite{voisin}.
		
		\begin{theorem}\label{thm: hodge}
		Let $X$ be a smooth projective variety over $\mathbb{C}$ (or more generally a compact K\"{a}hler manifold). For every $n$, the cohomology group $\H^n_{\B}(X)$ carries a functorial pure Hodge structure of weight $n$, where the subspace
		$$\H^{p,q}(X) \subset \H^n_\B(X)_\CC$$
		appearing in the Hodge decomposition is the space of cohomology classes that can be represented by a smooth differential form of type $(p,q)$.
		\end{theorem}
		
		Here $\H^n_\B(X)_\CC$ is identified with the cohomology of the complex of smooth differential forms (with complex coefficients) on $X^{\mathrm{an}}$, and a smooth differential form of type $(p,q)$ is one that can be written in local holomorphic coordinates $z_i$ as a linear combination of forms $f\,\d z_{i_1}\wedge \cdots \wedge \d z_{i_p}\wedge \d \overline{z}_{j_1}\wedge\cdots \wedge \d \overline{z}_{j_q}$ with $f$ a smooth function.
		
		Clearly, since $\H^2_\B(\mathbb{P}^1_\CC)$ has dimension $1$, we have the equality of pure Hodge structures
		$$\H^2_\B(\mathbb{P}^1_\CC)= \QQ(-1).$$
		More generally, for a connected smooth projective variety $X$ of dimension $n$ we have, by Poincar\'{e} duality:
		$$\H^{2n}_\B(X)=\QQ(-n).$$
		
	\subsection{The Hodge filtration}
		
		For a pure Hodge structure $H$ of weight $w$ we consider the decreasing filtration $\F$ on $H_\CC$, called the \emph{Hodge filtration}, defined by
		$$\F^pH_\CC := \bigoplus_{\substack{r+s=w\\ r\geq p}} H^{r,s}.$$
		One can recover the Hodge decomposition from the Hodge filtration by the formula
		$$H^{p,q} = \F^pH_\CC \cap \overline{\F^qH_\CC},$$
		and a pair $(H,\F)$ defines a pure Hodge structure of weight $w$ if and only if
		$$H_\CC = \F^pH_\CC \oplus \overline{\F^{w-p+1}H_\CC} \quad \mbox{ for all } p.$$
		
		The Hodge filtration is arguably better behaved and more fundamental than the Hodge decomposition. Indeed, it has a purely algebraic interpretation: via the isomorphism $\H^n_\B(X)_\CC\simeq \H^n_\dR(X)$,  the Hodge filtration on $\H^n_\dR(X)$ is induced by the ``stupid filtration''
		$$\F^p\Omega^\bullet_{X/\CC} := \Omega^{\bullet\geq p}_{X/\CC}$$
		on the complex of algebraic differential forms. As a consequence, it varies algebraically (or holomorphically) in families: if $\pi\colon X\to S$ is a smooth projective morphism of complex varieties over a smooth variety $S$, then the cohomology groups of the fibers of $\pi$ have pure Hodge structures by Theorem \ref{thm: hodge}, and the $\H^{p,q}$ of these fibers do not form an algebraic sub-bundle of $\mathcal{H}^n_\dR(X/S)$ (or a holomorphic sub-bundle of $\mathcal{H}^n_\dR(X/S)^{\mathrm{an}}$) in general, but the $\F^p$ do. These sub-bundles are not preserved by the Gauss--Manin connection but we have \emph{Griffiths transversality}:
		$$\nabla(\F^p\mathcal{H}^n_\dR(X/S)) \subset \F^{p-1}\mathcal{H}^n_\dR(X/S)\otimes_{\mathcal{O}_S} \Omega^1_{S/\CC}.$$
		This constraint motivates the definition of a \emph{variation of pure Hodge structure}, which we will not discuss here.

	\subsection{Mixed Hodge structures}
	
		The extension of classical Hodge theory to the cohomology of all complex varieties was achieved by Deligne \cite{delignehodge1,delignehodge2,delignehodge3}.
		
		\begin{definition}
		A \emph{mixed Hodge structure} is the datum of a finite dimensional vector space $H$ over $\QQ$ and
		\begin{enumerate}[$\bullet$]
		\item an increasing filtration $\W$ on $H$, called the \emph{weight filtration}
		\item a decreasing filtration $\F$ on $H_\CC$, called the \emph{Hodge filtration}
		\end{enumerate}
		such that for each integer $w$, the filtration induced by $\F$ on $\gr_w^\W H:=\W_wH/\W_{w-1}H$ defines a pure Hodge structure of weight $w$. A \emph{morphism of mixed Hodge structures} is a $\QQ$-linear map which is compatible with the weight filtrations and whose complexification is compatible with the Hodge filtrations.
		\end{definition}
		
		The category $\mathsf{MHS}$ of mixed Hodge structures is, surprisingly enough,  an abelian category. This comes from the fact that every morphism $f\colon H\to H'$ of mixed Hodge structures is \emph{strictly} compatible with the weight and Hodge filtrations, i.e., satisfies $f(\W_nH)=f(H)\cap \W_nH'$ and $f(\F^nH_\CC)=f(H_\CC)\cap \F^nH'_\CC$ for all $n$. The functors $H\mapsto \gr_n^\W H$ (resp. $H\mapsto \gr^n_\F H_\CC$) are exact functors from $\mathsf{MHS}$ to the category of finite dimensional vector spaces over $\QQ$ (resp. over $\CC$).
		
		\begin{remark}
		It is important to note that a mixed Hodge structure $(H,\W,\F)$ contains more information than just the collection of pure Hodge structures $\gr_w^\W H$ for all $w$. Indeed, the Hodge filtrations of those pure Hodge structures are not independent, but come from a \emph{single filtration} defined on $H$.
		\end{remark}
		
		The importance of mixed Hodge structures in the study of the topology of algebraic varieties comes from the following landmark result of Deligne.
		
		\begin{theorem}
		Let $X$ be an algebraic variety over $\CC$. For every $n$, the cohomology group $\H^n_\B(X)$ carries a functorial mixed Hodge structure, which is that of Theorem \ref{thm: hodge} in the smooth projective case.
		\end{theorem}
		
		More generally there is a functorial mixed Hodge structure on the relative cohomology groups $\H^n_\B(X,Y)$ for any pair $(X,Y)$ of complex varieties.
		
		\begin{example}
		Let $C$ be a smooth projective complex curve and let $a, b\in C$ be two distinct points. We have a short exact sequence
		\begin{equation}\label{eq: an extension MHS curves}
		0\To \H^1_\B(C)\To \H^1_\B(C\setminus \{a,b\}) \stackrel{r}{\To} \QQ(-1)\To 0
		\end{equation}
		in the category $\mathsf{MHS}$, where the map $r$ is dual to a small positively oriented cycle around $a$ on $C^{\mathrm{an}}$ and corresponds in de Rham cohomology to the residue of a $1$-form at $a$. One can also view $r$ as the connecting morphism in the Mayer--Vietoris long exact sequence for the covering of $C$ by the two opens $C\setminus \{a\}$ and $C\setminus \{b\}$, which explains why the target of $r$ is $\QQ(-1)=\H^2_\B(C)$. Therefore, the mixed Hodge structure on $H:=\H^1_\B(C\setminus \{a,b\})$ has two non-trivial weight-graded pieces, $\gr_1^\W H \simeq \H^1_\B(C)$ and $\gr_2^\W H  \simeq \QQ(-1)$. Even though the weight-graded quotients do not depend on $a$ and $b$, the mixed Hodge structure on $H$, i.e., the extension datum, does. In particular, $H$ is generally not isomorphic to the direct sum of $\H^1_\B(C)$ and $\QQ(-1)$ in $\mathsf{MHS}$. More concretely, the relevant extension group is
		$$\mathrm{Ext}^1_{\mathsf{MHS}}(\QQ(-1),\H^1_\B(C)) \simeq \mathrm{Pic}^0(C)_\QQ$$
		and the class of \eqref{eq: an extension MHS curves} corresponds to the difference $(a)-(b)\in\mathrm{Pic}^0(C)_\QQ$.
		\end{example}
	
		There is an obvious notion of \emph{tensor product} for mixed Hodge structures, which makes $\mathsf{MHS}$ into a neutral $\QQ$-linear tannakian category for which the forgetful functor $(H,\W,\F)\mapsto H$ is a fiber functor.
		
	\subsection{Mixed Hodge--Tate structures}
		
		\begin{definition}
		A \emph{mixed Hodge--Tate structure} is a mixed Hodge structure which is an iterated extension of the pure Hodge structures $\QQ(-n)$, for $n\in\mathbb{Z}$, or in other words one for which $\gr_w^\W$ vanishes for $w$ odd and is isomorphic to a direct sum of $\QQ(-n)$ for $w=2n$ even.
		\end{definition}
		
		An equivalent definition is as follows: a mixed Hodge--Tate structure is the datum of a finite dimensional vector space $H$ over $\QQ$ along with an increasing filtration $\W$ on $H$ and a decreasing filtration $\F$ on $H_\CC$ which satisfy
		$$\W_{2n}H_\CC = \W_{2(n-1)}H_\CC \oplus (\W_{2n}H_\CC\cap \F^nH_\CC) \quad \mbox{ for all } n.$$
		The Hodge filtration therefore splits the weight filtration over $\CC$. In other words, if we let $\gr^\W H:=\bigoplus_k\gr^\W_k H$, it induces an isomorphism
		\begin{equation}\label{eq: comp for MHTS}
		(\gr^\W H)_\CC \stackrel{\sim}{\To} H_\CC
		\end{equation}
		which sends $\bigoplus_{k\leq n}\gr_{2k}^\W H$ (resp. $\bigoplus_{k\geq n}\gr_{2k}^\W H$) to $\W_{2n}H$ (resp. to $\F^nH_\CC$) for all $n$.
		
		One conventionally renormalizes \eqref{eq: comp for MHTS} by multiplying it by $(2\pi\i)^k$ on $\gr_{2k}^\W H$ and considers its matrix in bases induced by a basis of $H$ compatible with the weight filtration. Such a matrix is block upper-triangular with diagonal blocks given by $(2\pi\i)^k$ times an identity matrix, where $k$ increases as one moves along the diagonal. It is called a \emph{period matrix} of the mixed Hodge--Tate structure, and completely determines its isomorphism class. For the sake of illustration, the following is a period matrix of a mixed Hodge--Tate structure $H$ for which $\gr^\W H\simeq \QQ(0)\oplus \QQ(-1)\oplus \QQ(-3)^{\oplus 2}$:
		
		$$\setlength{\arraycolsep}{3pt}\def\arraystretch{1.5}
		\left(\begin{matrix}  1 & * & * & * \\ 0 & 2\pi\i & * & * \\ 0 & 0 & (2\pi\i)^3 & 0 \\ 0 & 0 & 0 & (2\pi \i)^3\end{matrix}\right).$$
		
		\begin{remark}\label{rem: hodge splits weight}
		Such a matrix can be viewed as a period matrix in the sense of \S\ref{par: betti de rham} when the mixed Hodge--Tate structure $H$ comes from geometry, i.e., for instance, from the cohomology of a pair of algebraic varieties defined over $\QQ$. Indeed, in this case $H=M_\B$ has a de Rham counterpart $M_\dR$ which is a vector space over $\QQ$, equipped with a weight filtration $\W$ and a Hodge filtration $\F$, in such a way that the comparison isomorphism 
		$$\comp\colon M_\dR\otimes_\QQ\CC\stackrel{\sim}{\To} M_\B\otimes_\QQ\CC = H_\CC$$
		is compatible with the weight and Hodge filtrations. We insist on the fact that $\W$ and $\F$ are both defined over $\QQ$ on $M_\dR$, and therefore we get a canonical splitting
		$$\gr^\W M_\dR\stackrel{\sim}{\To} M_\dR$$
		which fits in the following commutative diagram:
		$$\xymatrixcolsep{1cm}\xymatrixrowsep{1cm}\diagram{
		(\gr^\W M_\dR)\otimes_\QQ\CC \ar[r] \ar[d]_-{\gr^\W(\comp)} & M_\dR \otimes_\QQ\CC \ar[d]^-{\comp} \\
		(\gr^\W H)_\CC \ar[r] & H_\CC
		}$$
		By Example \ref{ex: pi}, $\gr^\W(\comp)$ is given by multiplication by $(2\pi\i)^k$ in weight $2k$. Therefore the matrix of \eqref{eq: comp for MHTS} in a basis compatible with the weight filtration, suitably renormalized by powers of $2\pi\i$, is indeed the matrix of a comparison isomorphism between de Rham cohomology and Betti cohomology, and therefore deserves the name ``period matrix''. Note that this is very special to mixed Hodge--Tate structures, and in general a mixed Hodge structure does not contain enough information to recover a period matrix.
		\end{remark}
		
		The mixed Hodge--Tate structures of dimension $1$ are the pure Hodge--Tate structures $\QQ(-n)$, for $n\in\mathbb{Z}$. To classify the mixed Hodge--Tate structures of dimension $2$ it is enough to understand the extensions of $\QQ(-n)$ by $\QQ(0)$ for some $n\in\mathbb{Z}$, because the pure Hodge--Tate structures $\QQ(-i)$ are $\otimes$-invertible. One sees that we have
		\begin{equation}\label{eq: ext in MHTS}
		\operatorname{Ext}^1_{\mathsf{MHS}}(\QQ(-n),\QQ(0)) \simeq \begin{cases} \CC/(2\pi\i)^n\QQ & \mbox{ if } n\geq 1; \\ 0 & \mbox{ if } n\leq 0. \end{cases}
		\end{equation}
		For $n\geq 1$, the extension class corresponding to $\alpha\in\CC/(2\pi\i)^n\QQ$ has period matrix
		$$\setlength{\arraycolsep}{3pt}\def\arraystretch{1.3}
		\left(\begin{matrix}  1 & \alpha \\ 0 & (2\pi\i)^n\end{matrix}\right).$$
	
\section{Motives}

		We sketch the construction of several categories of motives. More details and references on the foundations of the theory of motives can be found in \cite{andrebook}.
	
		\subsection{The idea of motives}
		
			Grothendieck's idea of motives is that of a \emph{universal cohomology theory} for algebraic varieties. It involves an abelian $\QQ$-linear category $\Mot(F)$ whose objects are called \emph{motives over $F$} (with coefficients in $\QQ$), and functors
		\begin{equation}\label{eq: functor Var Mot}
		\M^n\colon (\mathsf{Var}/F)^{\operatorname{op}} \to \mathsf{Mot}(F) \; , \; X\mapsto \M^n(X)
		\end{equation}
		where $\mathsf{Var}/F$ denotes the category of varieties over $F$, and $\M^n(X)$ is called the \emph{motive of} $X$ in degree $n$. 
		
		One should think of $\M^n(X)$ as a universal object which controls the different cohomology groups $\H^n_?(X)$ for $?\in\{\dR,\B,\ldots\}$. More precisely, a cohomology theory, viewed as a functor $\H^\bullet_?\colon(\mathsf{Var}/F)^{\mathrm{op}} \to \mathsf{Vect}_k$ for some field $k$, should (under certain conditions) factor through \eqref{eq: functor Var Mot} via a \emph{realization functor}
		$$\omega_?\colon\mathsf{Mot}(F)\To \mathsf{Vect}_k.$$
		Furthermore, for a field $F\subset \CC$, the Betti realization functor $\omega_\B$
		should lift to a \emph{Hodge realization functor}
		$$\mathsf{Mot}(F)\To \mathsf{MHS}.$$
		This is sometimes called an \emph{enriched realization functor}.
				
		\subsection{Grothendieck's category of Chow motives}
		
			Grothendieck's construction of a category of motives (first written down by Manin \cite{manincorrespondences}) is only concerned with smooth projective varieties, and is based on algebraic cycles. The reason is that if $X$ and $Y$ are smooth projective varieties then an algebraic cycle $Z$ of codimension $r$ in $X\times Y$ (also known as a \emph{correspondence}) has a class $[Z]\in \H_?^{2r}(X\times Y)$ in any cohomology theory, and therefore induces a linear map
			\begin{equation}\label{eq: linear map from correspondence}
			\diagram{
			\H_?^\bullet(X) \ar[r]^-{(\operatorname{pr}_X)^*} & \H_?^\bullet(X\times Y)  \ar[r]^-{\cdot [Z]} & \H_?^{\bullet+2r}(X\times Y)\ar[r]^-{(\operatorname{pr}_Y)_*} & \H_?^{\bullet+2r-2\dim(X)}(Y),
			}
			\end{equation}
			where $\operatorname{pr}_X\colon X\times Y\to X$ and $\operatorname{pr}_Y\colon X\times Y\to Y$ denote the projections. Grothendieck's intuition is that those linear maps are the only ones that are common to all cohomology theories.
			 
			 The category of Chow motives with rational coefficients, denoted by $\mathsf{CHM}(F)$, is a symmetric monoidal $\QQ$-linear category equipped with a symmetric monoidal functor
			 $$\M \colon (\mathsf{SmProjVar}/F)^{\operatorname{op}} \To \mathsf{CHM}(F) \; , \; X\mapsto \M(X),$$
			 where the monoidal structure on the category $\mathsf{SmProjVar}/F$ of smooth projective varieties over $F$ is given by the product of varieties. Morphisms between Chow motives of varieties are related to algebraic cycles via the formula
			 \begin{equation}\label{eq: Hom in CHM}
			 \operatorname{Hom}_{\mathsf{CHM}(F)}(\M(X),\M(Y)) = \operatorname{CH}^{\dim(X)}(X\times Y)_\QQ,
			 \end{equation}
			 where $\operatorname{CH}^r$ denote the Chow groups of algebraic cycles of codimension $r$ modulo rational equivalence.  The composition of morphisms is given by the formula
			 $$Z_{23}\circ Z_{12} = (\operatorname{pr}_{13})_*\left((\operatorname{pr}_{12})^*(Z_{12})\cdot (\operatorname{pr}_{23})^*Z_{23}\right)$$
			 for algebraic cycles $Z_{ij}$ in $X_i\times X_j$, where $\operatorname{pr}_{ij}:X_1\times X_2\times X_3\to X_i\times X_j$ denote the projections.
			 
			 \begin{remark}\label{rem: CHM pseudo abelian}
			 Not all Chow motives are of the form $\M(X)$. Indeed, the category $\mathsf{CHM}(F)$ is \emph{pseudo-abelian} by design, and hence any projector $p$ of $\M(X)$ gives rise to a splitting $\M(X)=\operatorname{ker}(p)\oplus \operatorname{Im}(p)$ in $\mathsf{CHM}(F)$.
			 \end{remark}
			 
			 Every cohomology theory $\H^\bullet_?$ gives rise to a \emph{realization functor} on $\mathsf{CHM}(F)$ which sends $\M(X)$ to the direct sum of all the $\H^n_?(X)$ and an algebraic cycle $Z$ of codimension $\dim(X)$ in $X\times Y$ to \eqref{eq: linear map from correspondence}. For a field $F\subset \CC$ there is a Hodge realization functor
			 $$\mathsf{CHM}(F)\To \mathsf{MHS}$$
			 whose image lands in the category of direct sums of pure Hodge structures.
			 
			 \begin{remark}
			 One of the drawbacks of this construction is that it is not known whether a Chow motive $\M(X)$ splits as a direct sum of objects $\M^n(X)$ which lift the individual cohomology groups $\H_?^n(X)$, i.e., whether the ``K\"{u}nneth projectors'' $\H^\bullet_?(X)\twoheadrightarrow \H^n_?(X) \hookrightarrow \H^\bullet_?(X)$ are induced by a correspondence as in \eqref{eq: linear map from correspondence}. This is one of Grothendieck's \emph{standard conjectures on algebraic cycles} \cite{grothendieckstandard, kleimanstandard}.
			 \end{remark}

			 The unit of the tensor product of $\mathsf{CHM}(F)$ is $\QQ(0):=\M(\operatorname{Spec}F)$. More generally, we have for every integer $n$ a $\otimes$-invertible object $\QQ(-n)$ of $\mathsf{CHM}(F)$, and these objects satisfy $\QQ(-i)\otimes \QQ(-j)\simeq \QQ(-(i+j))$. One should think of $\QQ(-1)$ as playing the role of $\H^2(\PP^1)$, in the sense that we have a direct sum decomposition $\M(\PP^1_F) = \QQ(0)\oplus \QQ(-1)$, induced by any point $x$ of $\mathbb{P}^1_F$ as in Remark \ref{rem: CHM pseudo abelian} (the algebraic cycle $\mathbb{P}^1_F\times\{x\}\subset \mathbb{P}^1_F\times\mathbb{P}^1_F$ is a projector of $\M(\mathbb{P}^1_F)$). The tensor product 
			 $$M(-n):=M\otimes \QQ(-n)$$
			 is called a \emph{Tate twist} of $M$. Generalizing \eqref{eq: Hom in CHM} we have the following formula for morphisms between Tate twists of Chow motives of varieties:
			 $$\operatorname{Hom}_{\mathsf{CHM}(F)}(\M(X)(-m),\M(Y)(-n)) = \operatorname{CH}^{\dim(X)+m-n}(X\times Y)_\QQ.$$
			 It is justified by the fact that in general the composition \eqref{eq: linear map from correspondence} for $?=\B$ is a morphism of (pure) Hodge structures only up to the Tate twist $(-(\dim(X)-r))$ on the target.
		
		\subsection{Voevodsky's triangulated category of mixed motives}
		
		There does not seem to be any obvious way to generalize Grothendieck's construction of the category of Chow motives to incorporate varieties that are neither projective nor smooth. The first obvious obstruction is that the pushforward map $(\operatorname{pr}_Y)_*$ in \eqref{eq: linear map from correspondence} does not make sense in general.
		
		It was suggested by Beilinson and Deligne independently to ground the category of mixed motives on the \emph{complexes that compute cohomology groups} rather than on the cohomology groups themselves, following a two-step program:
		\begin{enumerate}[1)]
		\item Define a \emph{triangulated category} $\DM(F)$ which should play the role of the derived category of the abelian category of motives. This was done by Voevodsky \cite{voevodskytriangulated} and indepently by Hanamura \cite{hanamuramixed} and Levine \cite{levinemixed}. The category $\DM(F)$ is symmetric monoidal and equipped with a symmetric monoidal functor
		\begin{equation}\label{eq: functor from Var to DM}
		\M\colon \mathsf{Var}/F\To \DM(F) \; , \; X\mapsto \M(X).
		\end{equation}
		In the triangulated context, realization functors take values in derived categories of vector spaces $\mathsf{D}(\mathsf{Vect}_k)$. We warn the reader that \eqref{eq: functor from Var to DM} is naturally a \emph{covariant} functor, and therefore $\M(X)\in \DM(F)$ should be viewed as a motivic lift of a complex which computes the \emph{homology} of $X$.
		\item Extract the abelian category of motives using the Beilinson--Bernstein--Deligne--Gabber formalism of t-structures \cite{bbdg}. 	At present, this second step can only be made to work for certain subcategories of $\DM(F)$, as we will explain in the next section.
		\end{enumerate}
		
		Here is a rough idea of the steps involved in one possible definition\footnote{This is not Voevodsky's original definition, but gives an equivalent category. The category that we define here appears in the literature under the name \emph{category of \'{e}tale motives (without transfers)} and is usually denoted by $\mathsf{DA}_{\mathrm{\textnormal{\'{e}t}}}(F)$.} of $\DM(F)$:
		\begin{enumerate}[A)]
		\item $\,$We consider the category of presheaves of vector spaces over $\QQ$ on the category $\mathsf{SmVar}/F$ of smooth varieties defined over $F$, together with its (linearized) Yoneda embedding:
		\begin{equation}\label{eq: SmVar to PreSh}
		\mathsf{SmVar}/F \To \mathsf{PreSh}(\mathsf{SmVar}/F,\QQ) \; , \; X\mapsto \QQ\,\mathrm{Hom}_{\mathsf{SmVar}/F}(-,X) .
		\end{equation}
		We will be interested in those presheaves which are sheaves for the \'{e}tale topology on $\mathsf{SmVar}/F$, and in the sheafification functor
		\begin{equation}\label{eq: PreSh to Sh}
		\mathsf{PreSh}(\mathsf{SmVar}/F,\QQ) \To \mathsf{Sh}(\mathsf{SmVar}/F,\QQ).
		\end{equation}
		We denote the composition of \eqref{eq: SmVar to PreSh} and \eqref{eq: PreSh to Sh} by
		\begin{equation}\label{eq: yoneda embedding}
		\mathsf{SmVar}/F \To \mathsf{Sh}(\mathsf{SmVar}/F,\QQ) \; , \; X\mapsto \QQ(X).
		\end{equation}
		Finally, since we want to work with complexes, we consider the derived category of the abelian category $\mathsf{Sh}(\mathsf{SmVar}/F,\QQ)$.
		\item $\,$The second step is called \emph{$\mathbb{A}^1$-localization}. For a smooth variety $X$ over $F$, the projection map $X\times \mathbb{A}^1_F\to X$ induces an isomorphism in all cohomology theories, and we would like that to be true at the level of motives. We therefore consider the Verdier quotient of $\mathsf{D}(\mathsf{Sh}(\mathsf{SmVar}/F,\QQ))$ which forces the complexes
		$$0\to \QQ(X\times \AA^1_F)\To \QQ(X)\to 0$$
		to be zero. It is denoted by  $\DM^{\mathrm{eff}}(F)$ and called the category of \emph{effective Voevodsky motives} over $F$ (with rational coefficients). It is a $\QQ$-linear triangulated category and has a symmetric monoidal structure induced by the tensor product of vector spaces.
		\item $\,$Let $\QQ(1)$ be the object of $\mathsf{DM}^{\mathrm{eff}}(F)$ given by the complex $ \QQ(\PP^1_F)\To \QQ(\operatorname{Spec} F)$ induced by the structure morphism $\PP^1_F\to \operatorname{Spec}F$, where $\QQ(\PP^1_F)$ is in degree $-2$.  This should be thought of as playing the role of $\H_2(\PP^1_F)$. We define $\DM(F)$ to be the category obtained from $\DM^{\mathrm{eff}}(F)$ by formally adding a tensor inverse to $\QQ(1)$, denoted by $\QQ(-1)$.
		\item $\,$ Every smooth variety $X$ over $F$ gives rise to the object $\M(X)\in\DM(F)$ via \eqref{eq: yoneda embedding}, hence a functor $\mathsf{SmVar}/F\to \DM(F)$. With some more work, one can extend it to all varieties over $F$ and obtain \eqref{eq: functor from Var to DM}.
		\end{enumerate}
		
		The spaces of morphisms in $\DM(F)$ are related to algebraic cycles via a generalization of Chow groups discovered by Bloch and called \emph{higher Chow groups} \cite{blochhigher}, denoted by $\operatorname{CH}^r(X,n)$, where $\operatorname{CH}^r(X,0)=\operatorname{CH}^r(X)$ are the classical Chow groups. 
		We have an isomorphism:
		\begin{equation}\label{eq: Hom in DM}
		\operatorname{Hom}_{\DM(F)}(\M(X),\QQ(n)[i]) \simeq \operatorname{CH}^n(X,2n-i)_\QQ.
		\end{equation}
		
		\begin{remark}
		By looking at \eqref{eq: Hom in DM} for $i=2n$ one sees that the opposite category of $\mathsf{CHM}(F)$ is a full subcategory of $\DM(F)$. 
		\end{remark}		
		
		\begin{remark}
		The category of Chow motives is not abelian in general, and therefore is not the sought-for abelian category of pure motives, whose derived category is expected to embed in $\DM(F)$. Such a role is conjecturally played by a variant of Chow motives called \emph{numerical motives}, which is not based on Chow groups but rather on their quotients by numerical equivalence. We refer the reader to \cite[Chapitre 21]{andrebook} for more details on this matter.
		\end{remark}
	
		Voevodsky's triangulated category $\DM(F)$ is endowed with realization functors taking values in derived categories of vector spaces:
		$$\omega_?:\DM(F)\To \mathsf{D}(\mathsf{Vect}_k).$$
		One should think of $\omega_?(\QQ(X))$ as a complex which computes the \emph{homology} groups $\H_\bullet^?(X)$. For instance, if $F$ has an embedding $\sigma\colon F\to \CC$ then we have a Betti realization functor
		\begin{equation}\label{eq: Betti realization DM}
		\omega_{\B,\sigma}:\DM(F)\To \mathsf{D}(\mathsf{Vect}_\QQ)
		\end{equation}
		which is such that for every variety $X$ defined over $F$,
		$$\omega_{\B,\sigma}(\QQ(X)) = \operatorname{C}_\bullet^{\operatorname{sing}}(X_\sigma^{\mathrm{an}};\QQ)$$
		is the complex of rational singular chains on $X_\sigma^{\mathrm{an}}$.
		
\section{Second intermezzo: algebraic $\K$-theory}

	\subsection{Basics}
	
		Quillen \cite{quillenK} associates to every ring $R$ a sequence of abelian groups $\K_i(R)$, for $i\geq 0$, which are ``homotopical'' invariants of the category of finitely generated projective $R$-modules. We will here restrict to commutative rings $R$. The first three $\K$-groups, $\K_0$, $\K_1$, and $\K_2$, predate Quillen's work and are relatively easy to define (not so much to compute):
		\begin{enumerate}[$\bullet$]
		\item $\K_0(R)$ is the Grothendieck group of the category of finitely generated projective $R$-modules, i.e., the abelian group generated by isomorphism classes $[P]$ of finitely generated projective $R$-modules $P$, modulo the relations
		$$[P]=[P']+[P'']$$ 
		given by short exact sequences
		$$0\longrightarrow P'\longrightarrow P\longrightarrow P''\longrightarrow 0.$$ 
		The morphism $\ZZ\to \K_0(R) \, , \, 1\mapsto[R]$ is injective if $R$ is a non-zero commutative ring, and an isomorphism if all finitely generated projective $R$-modules are free, e.g., if $R$ is a field.\smallskip
		\item $\K_1(R)$ is the abelianization of the infinite general linear group $\operatorname{GL}(R)$, defined as the inductive limit of the groups $\operatorname{GL}_N(R)$ where $\operatorname{GL}_N(R)\hookrightarrow \operatorname{GL}_{N+1}(R)$ is given by $g\mapsto \left(\begin{matrix}g& 0 \\ 0 & 1\end{matrix}\right)$. There is a surjective map $K_1(R)\to R^\times$ given by the determinant of matrices. It is an isomorphism if $R$ is a field by Gaussian elimination.
		\item $\K_2(R)$ was defined by Milnor \cite{milnorintro}, and shortly after Matsumoto gave a presentation of $\K_2$ of a field by generators and relations \cite{matsumotoK2}.
		\end{enumerate}
	
		We refer the reader to \cite{rosenberg} for more information on the definition of higher $\K$-groups. They are notoriously difficult to compute, even in the case of a field. Their rationalized versions
		$$\K_i(R)_\QQ := \K_i(R)\otimes_\ZZ\QQ$$
		are more tractable objects because they can be computed in terms of group homology thanks to the Milnor--Moore theorem \cite{milnormoore}:
		$$\K_\bullet(R)_\QQ \simeq \operatorname{Prim}(\H_\bullet(\operatorname{GL}(R),\QQ)).$$
		(Here  $\H_\bullet(\operatorname{GL}(R),\QQ)$ has the structure of a graded Hopf algebra and $\operatorname{Prim}(\cdot)$ denotes its space of primitive elements, i.e., solutions of $\Delta(x)=1\otimes x+x\otimes 1$.)

	\subsection{The Dirichlet regulator and the analytic class number formula}
	
		Let $F$ be a number field and let $\mathcal{O}_F$ be its ring of integers. We have isomorphisms
		$$\K_0(\mathcal{O}_F)\simeq \ZZ\oplus \operatorname{Cl}(F),$$
		where $\operatorname{Cl}(F)$ is the ideal class group (a finite abelian group), and 
		$$\K_1(\mathcal{O}_F)\simeq \mathcal{O}_F^\times$$
		by a theorem of Bass--Milnor--Serre \cite{bassmilnorserre}. The latter group was computed by Dirichlet in 1846. Let $r_1$ denote the number of real embeddings $F\hookrightarrow \RR$, that we label $\sigma_1,\ldots,\sigma_{r_1}$, and $r_2$ denote the number of pairs of complex conjugate non-real embeddings $F\hookrightarrow \CC$, that we label $\sigma_{r_1+1},\ldots,\sigma_{r_1+r_2},\overline{\sigma_{r_1+1}},\ldots,\overline{\sigma_{r_1+r_2}}$. We have $[F:\QQ]=r_1+2r_2$. The \emph{Dirichlet regulator} is the group morphism
		$$\rho:\mathcal{O}_F^\times\oplus \ZZ \longrightarrow \RR^{r_1+r_2}$$
		defined by
		$$\rho(x)=(\log|\sigma_1(x)|,\ldots ,\log|\sigma_{r_1}(x)|,\log|\sigma_{r_1+1}(x)|^2,\ldots,\log|\sigma_{r_1+r_2}(x)|^2)$$ for $x\in\mathcal{O}_F^\times$, and $\rho(k)=\frac{1}{r_1+r_2}(k,\ldots,k)$ for $k\in\mathbb{Z}$. \emph{Dirichlet's unit theorem} states that $\rho$ is injective modulo torsion and that its image is a lattice in $\RR^{r_1+r_2}$. In particular,
		$$\operatorname{rk} \mathcal{O}_F^\times = r_1+r_2-1.$$
		
		The \emph{regulator} $R_F$ is defined as the covolume of the image of $\rho$ in $\RR^{r_1+r_2}$. The \emph{analytic class number formula} expresses the residue of the Dedekind zeta function of $F$ at $s=1$ in terms of the regulator, the number $w_F$ of roots of unity in $F$, the discriminant $D_F$, and the class number $h_F=|\operatorname{Cl}(F)|$:
		$$\lim_{s\to 1}(s-1)\zeta_F(s) = \frac{2^{r_1+r_2}h_F}{w_F}\frac{\pi^{r_2}}{\sqrt{|D_F|}}R_F.$$
		In what follows we will not be concerned with rational prefactors and will therefore write the formula as:
		\begin{equation}\label{eq: class number formula explicit}
		\lim_{s\to 1}(s-1)\zeta_F(s) \sim_{\QQ^\times } \frac{\pi^{r_2}}{\sqrt{|D_F|}} \det\!\left(\log|\sigma_i(\varepsilon_j)|\right)_{1\leq i,j\leq r_1+r_2-1},
		\end{equation}
		where $(\varepsilon_1,\ldots,\varepsilon_{r_1+r_2-1})$ is a basis of $\mathcal{O}_F^\times$ modulo torsion, and $a\sim_{\QQ^\times} b$ means that $a\in\QQ^\times b$.
	
	\subsection{The Borel regulator and higher $\K$-theory of number fields}
	
		Borel \cite{borelstable, borelzeta} made crucial contributions to the computation of the algebraic $\K$-theory of number fields, in the spirit of Dirichlet's work and the analytic class number formula. For an integer $n\geq 2$ let us write 
		$$d_n = \begin{cases} r_1+r_2 & \mbox{ if } n \mbox{ is odd} \\ r_2 & \mbox{ if } n \mbox{ is even} \end{cases}$$
		and fix the following identification, induced by the real or imaginary part:
		\begin{equation}\label{eq: identification dn}
		\Bigg(\bigoplus_{\sigma\colon F\to \CC}\CC/(2\pi\i)^n\RR\Bigg)^+\simeq \RR^{d_n},
		\end{equation}
		where the symbol $+$ denotes the space of invariants for complex conjugation acting on each $\CC/(2\pi\i)^n\RR$ and on the set of embeddings $\sigma\colon F\to \CC$.
		
		Borel proved that all even rational $\K$-groups $\K_{2i}(F)_\QQ$, for $i\geq 1$,  are zero. In the case of odd $\K$-groups he defined for all $n\geq 2$ a class in $\H^{2n-1}(\operatorname{GL}(\CC),\CC/(2\pi\i)^n\RR)$, giving rise to a morphism 
		$$\K_{2n-1}(\CC)\to \CC/(2\pi\i)^n\RR,$$ 
		which is compatible with complex conjugation. By functoriality of $\K$-theory for all embeddings $\sigma\colon F\to \CC$, and taking \eqref{eq: identification dn} into account, we obtain a morphism
		\begin{equation}\label{eq: borel regulator}
		\rho_n:\K_{2n-1}(F)\longrightarrow\RR^{d_n},
		\end{equation}
		called the \emph{Borel regulator}, which is a higher version of the Dirichlet regulator\footnote{The Dirichlet regulator is related to $\K_1(\mathcal{O}_F)\simeq \mathcal{O}_F^\times$ and not $\K_1(F)\simeq F^\times$, which has infinite rank. However, for $i\geq 2$ we have $K_i(\mathcal{O}_F)_\QQ\simeq K_i(F)_\QQ$.}, as the following important theorem of Borel shows.
		
		\begin{theorem}\label{thm: borel}
		Let $n\geq 2$. The Borel regulator $\rho_n$ is injective modulo torsion and its image is a lattice in $\RR^{d_n}$. In particular 
		\begin{equation}\label{eq: rank K theory}
		\operatorname{rk}\K_{2n-1}(F)=d_n.
		\end{equation}
		Furthermore, the covolume $R_F^{(n)}$ of the image of $\rho_n$ in $\RR^{d_n}$ is related to the value at $s=n$ of the Dedekind zeta function of $F$ via the formula:
		$$\zeta_F(n) \sim_{\QQ^\times} \frac{\pi^{n([F:\QQ]-d_n)}}{\sqrt{|D_F|}}R_F^{(n)}.$$
		\end{theorem}
	
		In contrast with the Dirichlet regulator, the Borel regulators are not easy to express explicitly in terms of special functions, even for small $n$ and small number fields. As we will now see, a more concrete description of Borel regulators would follow from a deep understanding of categories of mixed Tate motives.

\section{Mixed Tate motives}

	\subsection{The abelian category of mixed Tate motives: definition}
	
		\begin{definition}\label{defi: DMT}
		Let $F$ be a field. The \emph{triangulated category of mixed Tate motives over $F$}, denoted by~$\DMT(F)$, is the triangulated subcategory of $\DM(F)$ generated by the pure Tate motives $\QQ(-n)$, for $n\in\ZZ$.
		\end{definition}
		
		In other words, an object of $\DMT(F)$ is an iterated extension (in the triangulated sense) of the objects $\QQ(-n)[r]$, for $n,r\in\ZZ$. The key to understanding mixed Tate motives therefore lies in the morphisms between shifts of pure Tate motives. Generalizing the Riemann--Roch isomorphism involving Chow groups and $\K_0$, Bloch \cite{blochhigher} proved a general relationship between higher Chow groups and algebraic $\K$-theory, which allows one to rewrite \eqref{eq: Hom in DM} for $X=\operatorname{Spec} F$ as follows:
		\begin{equation}\label{eq: Hom in DMT}
		\operatorname{Hom}_{\DMT(F)}(\QQ(-n),\QQ(0)[i]) \simeq \gr^n_\gamma \K_{2n-i}(F)_\QQ.
		\end{equation}
		Here $\gr_\gamma$ denotes the successive quotients for the $\gamma$-filtration defined by Soul\'{e} \cite{souleoperations}, a byproduct of the $\lambda$-structure in $\K$-theory (or, alternatively, of the Adams operations). The following vanishing conjecture, stated independently by Beilinson \cite{beilinsonhigher} and Soul\'{e} \cite{souleoperations}, plays a central role in the study of mixed Tate motives.
	
		\begin{conjecture}\label{conj: BS}
		For any field $F$ and any integer $n\geq 1$, 
		\begin{equation}\label{eq: vanishing BS}
		\gr^n_\gamma\K_{2n-i}(F)_\QQ = 0 \quad \mbox{ if }\; i\leq 0.
		\end{equation}
		\end{conjecture}
		
		This conjecture is known to hold for finite fields, whose higher rational $\K$-groups vanish by Quillen \cite{quillenfinitefields}, and, more importantly for us, for number fields by the work of Borel discussed in the previous section. Indeed, in this case we have, for all $n\geq 1$:
		\begin{equation}\label{eq: gamma filtration number field}
		\gr_\gamma^n\K_{2n-i}(F)_\QQ= 0 \quad \mbox{ if }\; i\neq 1.
		\end{equation}
		
		\begin{definition}
		Let $F$ be a field for which the Beilinson--Soul\'{e} vanishing \eqref{eq: vanishing BS} holds. The \emph{abelian category of mixed Tate motives over $F$}, denoted by~$\MT(F)$, is the full subcategory of $\DMT(F)$ consisting of iterated extensions (in the triangulated sense) of the pure Tate motives $\QQ(-n)$, for $n\in\ZZ$.
		\end{definition}
		
		The difference with Definition \ref{defi: DMT} is that shifts are not allowed. Levine proves \cite{levinetate} that $\MT(F)$ is an abelian category by recasting it as the heart of a natural t-structure on $\DMT(F)$. The Beilinson--Soul\'{e} vanishing \eqref{eq: vanishing BS} enters the proof of the vanishing of $\operatorname{Hom}_{\DMT(F)}(M,M')$ for $M\in \DMT(F)^{\leq 0}$ and $M'\in \DMT(F)^{\geq 1}$, which is one of the axioms of a t-structure (mimicking the fact that for objects $A$, $B$ of an abelian category and an integer $i<0$, every morphism $A\to B[i]$ in the derived category is zero).

	\subsection{The abelian category of mixed Tate motives: formal structure}

		The category $\MT(F)$ has the same formal properties as the category of mixed Hodge--Tate structures:
		\begin{enumerate}[1)]
		\item All the pure Tate objects $\QQ(-n)$ are in $\MT(F)$, and they are the only simple objects. We have
		$$\operatorname{Hom}_{\MT(F)}(\QQ(-n),\QQ(-n')) = \begin{cases} \QQ & \mbox{ if } n=n';\\ 0 & \mbox{ otherwise.}\end{cases}$$
		\item Every object $M\in \MT(F)$ is equipped with a finite \emph{weight filtration} $\W$ indexed by even integers, such that $\gr_{2n}^\W M$ is isomorphic to a direct sum of a finite number of copies of $\QQ(-n)$. The weight filtration is functorial.
		\item The tensor product of $\DM(F)$ makes $\MT(F)$ into a neutral tannakian category over $\QQ$, for which
		\begin{equation}\label{eq: fiber functor MT}
		\omega:\MT(F)\longrightarrow \mathsf{Vect}_\QQ \; , \; M\mapsto \bigoplus_{n\in\ZZ}\operatorname{Hom}_{\MT(F)}(\QQ(-n),\gr_{2n}^\W M)
		\end{equation}
		is a fiber functor.
		\end{enumerate}
		
		We also have at our disposal the classical de Rham and Betti fiber functors
		$$\omega_\dR \colon \MT(F)\To \mathsf{Vect}_F \quad \mbox{ and } \quad \omega_{\B,\sigma}\colon \MT(F)\To \mathsf{Vect}_\QQ,$$
		when they make sense. The ``canonical'' fiber functor \eqref{eq: fiber functor MT} is a rational structure on the de Rham fiber functor
		\begin{equation}\label{eq: omega dR vs omega canonical}
		\omega_\dR\simeq \omega\otimes_\QQ F,
	\end{equation}
	which essentially follows from the argument of Remark \ref{rem: hodge splits weight}.

	\subsection{The case of number fields, and the Hodge regulator}
		
		If $F$ is a number field then \eqref{eq: Hom in DMT} and \eqref{eq: gamma filtration number field} yield a complete description of the Ext groups between pure Tate motives:
		\begin{equation}\label{eq: ext in MT F number field}
		\operatorname{Ext}^1_{\MT(F)}(\QQ(-n),\QQ(0)) \simeq \begin{cases} \K_{2n-1}(F)_\QQ & \mbox{ if } n\geq 1; \\ 0 & \mbox{ if } n\leq 0.\end{cases}
		\end{equation}
		Furthermore, all higher Ext groups vanish. There are two different cases to consider:
		\begin{enumerate}[$\bullet$]
		\item For $n=1$ we get
		\begin{equation}\label{eq: ext MT F n=1}
		\operatorname{Ext}^1_{\MT(F)}(\QQ(-1),\QQ(0)) \simeq \K_1(F)_\QQ \simeq F^\times_\QQ,
		\end{equation}
		which has infinite dimension. Those extensions of $\QQ(-1)$ by $\QQ(0)$ are easy to describe: they are the \emph{Kummer motives}, which lift the objects of Example \ref{ex: log}, as we will see in the next paragraph.
		\item For $n\geq 2$, \eqref{eq: ext in MT F number field} has finite dimension, given by Borel's theorem \eqref{eq: rank K theory}. However, the corresponding extensions are \emph{difficult to describe explicitly} even for small number fields $F$. An explicit description of those extensions would lead to explicit formulas for Borel's regulator, in the following sense. 
		For every complex embedding $\sigma\colon F\to \CC$ there is a Hodge realization functor \cite{delignegoncharov}
		$$\MT(F)\longrightarrow \mathsf{MHTS},$$
		which sends the mixed Tate motive $\QQ(-n)$ to the mixed Hodge--Tate structure $\QQ(-n)$. In view of \eqref{eq: ext in MHTS} and \eqref{eq: ext in MT F number field}, the induced morphisms on $\operatorname{Ext}^1$ read
		$$\varpi_n^{(\sigma)}:\K_{2n-1}(F)_\QQ \To \CC/(2\pi\i)^n\QQ.$$
		 By summing them over all complex embeddings of $F$ and composing with the quotient map $\CC/(2\pi\i)^n\QQ\to \CC/(2\pi\i)^n\RR$, we obtain via \eqref{eq: identification dn} the \emph{Hodge regulator}
		 $$\varpi_n:\K_{2n-1}(F)_\QQ\To \RR^{d_n},$$
		 which is expected to be equal to the Borel regulator \eqref{eq: borel regulator}:
		 \begin{equation}\label{eq: equality of regulators}
		 \rho_n\stackrel{?}{=}\varpi_n.
		 \end{equation}
		\end{enumerate}
		
	\subsection{A first family of mixed Tate motives: Kummer motives}
	
		We construct motivic lifts of the objects of Example \ref{ex: log}. Let $F$ be a number field, let $x\in F^\times$, and consider the object $M$ of $\DM(F)$ defined by the complex
		$$0\to \QQ(\{1\})\oplus \QQ(\{x\}) \longrightarrow \QQ(\mathbb{A}^1_F\setminus \{0\})\to 0,$$
		where $\QQ(\mathbb{A}^1_F\setminus \{0\})$ is placed in degree $-1$ and the differential is induced by the closed immersions of $\{1\}$ and $\{x\}$ inside $\mathbb{A}^1_F\setminus\{0\}$. It plays the role\footnote{Strictly speaking, this is only true if $x\neq 1$.} of the relative \emph{homology} group $\H_1(\mathbb{A}^1\setminus\{0\},\{1,x\})$. The object of $\DM(F)$ defined by the subcomplex
		$$\QQ(\{1\})\longrightarrow \QQ(\mathbb{A}^1_F\setminus \{0\})$$
		is isomorphic to $\QQ(1)$ by a motivic lift of the Mayer--Vietoris argument that proves that $\H_1(\mathbb{A}^1\setminus\{0\})\simeq \H_2(\PP^1)$. Since the object $\QQ(\{x\})$ is the trivial object $\QQ(0)$, we see that $M$ sits in a distinguished triangle
		$$\QQ(1) \longrightarrow M \longrightarrow \QQ(0) \stackrel{+1}{\longrightarrow}$$
		in $\DM(F)$. This proves that $M$ is an object of the category $\MT(F)$. Its dual, which plays the role of the relative \emph{cohomology} group $\H^1(\mathbb{A}^1\setminus\{0\},\{1,x\})$ is called the \emph{Kummer motive} $\mathcal{K}_x$, and sits in a short exact sequence
		\begin{equation}\label{eq: short exact sequence Kummer motive}
		0\longrightarrow \QQ(0) \longrightarrow \mathcal{K}_x \longrightarrow \QQ(-1)\longrightarrow 0
		\end{equation}
		in $\MT(F)$. Under the isomorphism \eqref{eq: ext MT F n=1}, the extension class of $\mathcal{K}_x$ corresponds to the element $x\in F^\times_\QQ$, and we see that we completely understand the extensions of $\QQ(-1)$ by $\QQ(0)$ in $\MT(F)$. As seen in Example \ref{ex: log}, the period matrix of $\mathcal{K}_x$, relative to an embedding $\sigma\colon F\to \CC$, is \eqref{eq: period matrix log} with $z=\sigma(x)$, and the corresponding class in $\operatorname{Ext}^1_{\mathsf{MHTS}}(\QQ(-1),\QQ(0))\simeq \CC/2\pi\i\QQ$ is therefore $\log(\sigma(x))$.
		
	\subsection{What are the mixed Tate motives of rank $2$?}
	
		Let $F$ be a number field, let $n\geq 2$, and assume that we are given a family of extensions
		\begin{equation}\label{eq: some extension in MT F}
		0\To \QQ(0) \To E_j \To \QQ(-n)\To 0
		\end{equation}
		in $\MT(F)$, for $j=1,\ldots,d_n$, whose classes span the corresponding Ext group. Write a period matrix of $E_j$ relative to an embedding $\sigma\colon F\to \CC$ as
		$$\setlength{\arraycolsep}{3pt}\def\arraystretch{1.3}
		\left(\begin{matrix}  1 & \alpha_{j}^{(\sigma)} \\ 0 & (2\pi\i)^n\end{matrix}\right).$$
		Then, under the equality of regulators \eqref{eq: equality of regulators} we have, thanks to Theorem \ref{thm: borel}, an expression for the special value of the Dedekind zeta function of $F$:
		\begin{equation}\label{eq: general formula for zeta F n}
		\zeta_F(n)\sim_{\QQ^\times} \frac{\pi^{n([F:\QQ]-d_n)}}{\sqrt{|D_F|}}\det\left(\alpha_j^{(\sigma_i)}\right)_{1\leq i,j\leq d_n}.
		\end{equation}
		Such a formula would be a higher version of the analytic class number formula \eqref{eq: class number formula explicit}.
		
		Unfortunately, the extensions \eqref{eq: some extension in MT F} are hard to find ``in nature'', i.e., when looking at relative cohomology groups of pairs of varieties. We break the quest for these extensions into two subtasks:
		\begin{enumerate}[1)]
		\item Produce interesting families of mixed Tate motives over $F$ that are ``as rich as possible''. They will almost never be extensions \eqref{eq: some extension in MT F}. Section \ref{sec: families MT} will be devoted to a description of such families.
		\item Understand those examples well enough to create extensions \eqref{eq: some extension in MT F} out of them. The main tool for this is the tannakian formalism, that we will develop in sections \ref{sec: tannakian} and \ref{sec: tannakian examples}.
		\end{enumerate}
		
\section{Some families of mixed Tate motives}\label{sec: families MT}

	\subsection{How to construct mixed Tate motives?}
	
		The Voevodsky motive $\M(X)\in \DM(F)$ of a variety $X$ over $F$ almost never lives in the triangulated category $\DMT(F)$ of mixed Tate motives, for the same reason that the mixed Hodge structure on $\H^\bullet_\B(X)$ is almost never a mixed Hodge--Tate structure. However, the motive of affine $n$-space is $\M(\mathbb{A}^n_F)\simeq \M(\operatorname{Spec} F)=\QQ(0)$, and is therefore in $\DMT(F)$. One can use this fact to construct more mixed Tate motives, as follows.
		
		If $X$ is smooth and decomposes as the union $Z\cup U$ where $Z$ is a smooth subvariety of codimension $r$ and $U$ is its complement, then the \emph{localization triangle}
		\begin{equation}\label{eq: localization triangle}
		\M(Z)[-2r](-r)\To \M(X)\To \M(U)\stackrel{+1}{\To}
		\end{equation}
		shows that if two out of $\M(Z)$, $\M(X)$, $\M(U)$ are in $\DMT(F)$ then the third one is as well.  Hence, the motive $\M(X)$ of a variety $X$ obtained by ``cutting and pasting'' affine spaces is in $\DMT(F)$. Taking cohomology with respect to Levine's t-structure and dualizing therefore produces objects
		$$\M^n(X) \in \MT(F)$$
		which lift the cohomology groups $\H^n_?(X)$.
		Unfortunately, such mixed Tate motives are rarely very interesting in the sense that they tend to be \emph{pure} Tate motives in a lot of cases. For instance, for any union $L$ of hyperplanes in $\PP^N_F$, the motive $\M^n(\PP^N_F\setminus L)$ is a direct sum of the pure Tate motives $\QQ(-n)$, as can be seen by induction on the number of hyperplanes using the localization triangle \eqref{eq: localization triangle}. 
		
		The key to producing \emph{interesting} (i.e., non-pure) mixed Tate motives is to consider \emph{relative} cohomology groups. Let $(X,Y)$ be a pair defined over $F$, write $Y$ as the union of varieties $Y_i$ and assume that $X$, the $Y_i$, and all their multiple intersections are smooth. Consider the complex
		\begin{equation}\label{eq: complex relative cohomology DM}
		\M(X,Y) : \; \cdots \To \bigoplus_{i<j} \QQ(Y_i\cap Y_j) \To \bigoplus_i \QQ(Y_i) \To \QQ(X)\To 0,
		\end{equation}
		whose arrows are alternating sums of the maps induced by the natural closed immersions.  If all the objects $\M(X)$, $\M(Y_i)$, $\M(Y_i\cap Y_j)$, and so on, are in $\DMT(F)$, then $\M(X,Y)$ is as well. Taking cohomology for Levine's t-structure and dualizing produces objects
		$$\M^n(X,Y)\in \MT(F)$$
		which lift the relative cohomology groups $\H^n_?(X,Y)$. In the case $X=\AA^1_F\setminus\{0\}$ and $Y=\{1,x\}$, this is how we defined Kummer motives in the previous section.
		
	\subsection{Polylogarithm motives}
	
	For $x\in F\setminus \{0,1\}$ we define the $n$-th \emph{polylogarithm motive}	
	$$\mathcal{L}_x^{(n)} = \M^n(\mathbb{A}^n_F\setminus \{xx_1\cdots x_n=1\}, \{x_1(1-x_1)\cdots x_n(1-x_n)=0\}),$$
	which lifts \eqref{eq: relative cohomology polylog}. As explained in \cite{dupontfresan}, it is an object of $\MT(F)$ whose period matrix, relative to an embedding $\sigma\colon F\to \CC$, is \eqref{eq: period matrix polylog} with 
	$$z=\sigma(x).$$ 
	We have short exact sequences
	\begin{equation}\label{eq: polylog as extension}
	0\To \QQ(0)\To \mathcal{L}_x^{(n)} \To \operatorname{Sym}^{n-1}(\mathcal{K}_x)(-1)\To 0
	\end{equation}
	and
	 \begin{equation}\label{eq: polylog inductive}
	0\To \mathcal{L}_x^{(n-1)}\To \mathcal{L}_x^{(n)} \To \QQ(-n)\To 0,
	\end{equation}
	which reflect the block-triangular shape of \eqref{eq: period matrix polylog}.
	
	\subsection{Motives of bi-arrangements of hyperplanes, and Aomoto polylogarithms}
	
	Kummer motives are of the form
	$$\M^1(\PP^1_F\setminus \{0,\infty\},\{1,x\})$$
	and one can naturally construct more general mixed Tate motives using projective geometry in higher dimension. From the data $(L,M)$ of two unions of hyperplanes in projective space $\mathbb{P}^n_F$, one constructs mixed Tate motives
	\begin{equation}\label{eq: relative cohomology aomoto generic}
	\M^n(\PP^n_F\setminus L,M\setminus L\cap M).
	\end{equation}
	Their periods, relative to an embedding of $F$ in $\CC$, look like
	$$\int_\alpha\omega,$$
	where $\omega$ is an algebraic $n$-form on $\PP^n_F\setminus L$, and $\alpha$ is a topological $n$-cycle in $(\PP^n\setminus L)(\CC)$ whose boundary lies in $M(\CC)$.
	
\begin{figure}[ht]
	\centering 
	\subfigure{\begin{tikzpicture}
	\def\a{2};
	\draw[thick,VioletRed3] (0,-.3*\a) -- (0,2*\a);
	\draw[thick,VioletRed3] (-.25*\a,-.25*\a) -- (1.8*\a,1.8*\a);
	\draw[thick,VioletRed3] (-.3*\a,\a) -- (2*\a,\a);
	\draw[thick,VioletRed3,fill=VioletRed3!20!white] (0,0) -- (\a,\a) -- (0,\a) -- cycle;
	\draw[thick,Turquoise4] (-.3*\a,0) -- (2*\a,0);
	\draw[thick,Turquoise4] (1.3*\a,-.3*\a) -- (1.3*\a,2*\a);
	\node[circle,inner sep=2pt,draw=VioletRed3,fill=VioletRed3, thick] (root) at (0,0) {};
	\end{tikzpicture}}
	 \hspace{1.5cm} 
 	\subfigure{\begin{tikzpicture}
	\def\a{2};
	\draw[thick,VioletRed3] (0,.2*\a) -- (0,2*\a);
	\draw[thick,VioletRed3] (45:.2*\a) -- (1.8*\a,1.8*\a);
	\draw[thick,VioletRed3] (-.3*\a,\a) -- (2*\a,\a);
	\draw [VioletRed3,thick,domain=-30:120] plot ({.3*\a*cos(\x)}, {.3*\a*sin(\x)});
	\draw[thick,VioletRed3,fill=VioletRed3!20!white] (45:.3*\a) arc (45:90:.3*\a) -- (0,.3*\a) -- (0,\a) -- (\a,\a) -- cycle;
    (90:3mm) arc (90:0:3mm) arc (180:90:3mm) ;
	\draw[thick,Turquoise4] (.2*\a,0) -- (2*\a,0);
	\draw[thick,Turquoise4] (1.3*\a,-.3*\a) -- (1.3*\a,2*\a);
	\end{tikzpicture}}
	 \caption{In $\mathbb{A}^2$, the blow-up of the origin separates the boundary of the simplex $\{0\leq t_1\leq t_2\leq 1\}$ from the pole divisor $\{(1-zt_1)t_2=0\}$. By removing the strict transform of $\{t_2=0\}$ one recovers the geometry of Figure \ref{fig: polylog}.}	
	 \label{fig: blow up}
	\end{figure}
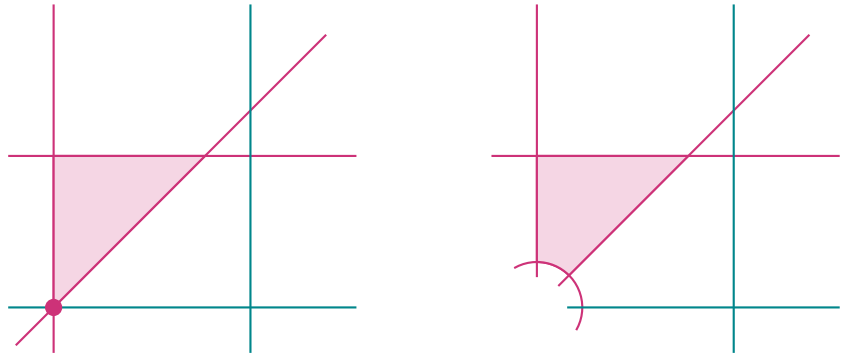

	This construction is not very useful when $L,M$ are not in general position because in this case there might not be many interesting integration domains $\alpha$. For instance, working in affine space instead of projective space for simplicity, consider the geometric situation of Figure \ref{fig: blow up}, with $L=\{(1-zt_1)t_2=0\}$ and $M=\{t_1(t_2-t_1)(1-t_2)=0\}$. The triangle $\{0\leq t_1\leq t_2\leq 1\}$ meets the line $\{t_2=0\}$ and therefore does not live in $(\PP^n\setminus L)(\CC)$. However, the integral
	\begin{equation}\label{eq: Li 2 simplicial}
	\iint_{0\leq t_1\leq t_2\leq 1}\frac{z\,\d t_1}{1-zt_1}\frac{\d t_2}{t_2}
	\end{equation}
	converges if $z\notin (1,\infty)$, and in fact equals $\Li_2(z)$ (expand $\frac{1}{1-z t_1}$ as a geometric series and integrate first in $t_1$ then in $t_2$).  In order to view it as a period of relative cohomoloy, we work in the blow-up $\pi:X\to \AA^2$ along the triple point $(0,0)$ as in Figure \ref{fig: blow up} below. There, the pull-back $\pi^*(\omega)$ only has poles along the strict transform $A=\widetilde{L}$, and the triangle $\{0\leq t_1\leq t_2\leq 1\}$ transforms into a quadrilateral with boundary along the union $B=\widetilde{M}\cup E$ where $E$ is the exceptional divisor. The motive
	$$\M^2(X\setminus A,B\setminus A\cap B)$$
	therefore has \eqref{eq: Li 2 simplicial} as a period.
	
	The generalization of this construction to all pairs $(L,M)$ leads to interesting mixed Tate motives which we call \emph{motives of bi-arrangements of hyperplanes} (some tools for computing them can be found in \cite{dupontbiarrangements}). They were considered for the first time by Beilinson \emph{et al.} in \cite{bvgs}, where their periods are called \emph{Aomoto polylogarithms} as an homage to Aomoto's work on related families of integrals and cohomology groups (see, e.g., \cite{aomotostructure}). Polylogarithm motives are special cases of this construction.

	\subsection{Iterated integrals and motivic fundamental groups}
	
	Let $X$ be a manifold, let $\gamma:[0,1]\to X$ be a smooth path, and let $\omega_1,\ldots,\omega_n$ be smooth $1$-forms on $X$. We define the \emph{iterated integral} of the formal word $\omega_1\cdots \omega_n$ along $\gamma$ as
	\begin{equation}\label{eq: iterated integral}
	\int_\gamma \omega_1\cdots \omega_n := \int_{0\leq t_1\leq\cdots \leq t_n\leq 1} f_1(t_1) \,\d t_1\cdots f_n(t_n) \,\d t_n
	\end{equation}
	where we have set $f_i(t)\,\d t:=\gamma^*(\omega_i)$.  The case $n=1$ is simply the integral of a $1$-form along a smooth path. If $X$ is the analytification of an algebraic variety and the $\omega_i$'s are algebraic, then \eqref{eq: iterated integral} is a period of relative cohomology, by a theorem of Beilinson \cite{goncharovmultiple, delignegoncharov}. We describe the corresponding Betti--de Rham comparison isomorphism in a special case.
	
	Let $S\subset \mathbb{P}^1(F)$ be a finite set of points, with $\infty\in S$, and let us consider the punctured projective line $X=\mathbb{P}^1_F\setminus S$. We fix two basepoints $a,b\in X(F)$. Deligne and Goncharov define in \cite{delignegoncharov} an ind-mixed Tate motive over $F$ (i.e., an object of the ind-category of $\MT(F)$), called the (algebra of functions on the) \emph{motivic fundamental group}, that we denote by
	$$M=\mathcal{O}(\pi_1^{\operatorname{mot}}(X, a,b)).$$
	\begin{enumerate}[$\bullet$]
	\item Its de Rham realization is
	$$M_\dR = F\langle \omega_s, s\in S\setminus\{\infty\}\rangle,$$
	the space of words in the letters
	$$\omega_s = \frac{\d x}{x-s} \quad (s\in S\setminus\{\infty\})$$
	which form a basis of $\H^1_\dR(X)$. The weight filtration remembers the length of the words.
	\item Its Betti realization relative to an embedding of $F$ in $\CC$ is the algebra of functions on the \emph{unipotent completion} of $\pi_1(X^{\mathrm{an}},a,b)$. To define it, let $I\subset \QQ[\pi_1(X^{\mathrm{an}},a)]$ denote the augmentation ideal of the group algebra of the topological fundamental group of $X^{\mathrm{an}}$ based at $a$, and consider its powers $I^{n+1}$, for $n\geq 0$. Any path $\delta$ from $a$ to $b$ in $X^{\mathrm{an}}$ induces a $\QQ$-linear isomorphism $\QQ[\pi_1(X^{\mathrm{an}},a)]\simeq \QQ[\pi_1(X^{\mathrm{an}},a,b)],\gamma\mapsto \gamma\delta$, and we can view the ideals $I^{n+1}$ as subspaces of $\QQ[\pi_1(X^{\mathrm{an}},a,b)]$, which are independent of the choice of $\delta$. We then have
	$$M_\B = \lim_n \left(\QQ[\pi_1(X^{\mathrm{an}},a,b)] / I^{n+1}\right)^\vee,$$
	the space of $\QQ$-valued functions on $\pi_1(X^{\mathrm{an}},a,b)$ which are zero on $I^{n+1}$ for some $n\geq 0$. The weight filtration remembers the index $n$.
	\item The Betti--de Rham comparison
	\begin{equation}\label{eq: comp for pi_1}
	M_\dR\otimes_F\CC\stackrel{\sim}{\To} M_\B\otimes_\QQ\CC
	\end{equation}
	is induced by iterated integration, i.e., by 
	$$\pi_1(X^{\mathrm{an}},a,b) \times F\langle \omega_s,s\in S\setminus\{\infty\} \rangle \To \CC \; , \; (\gamma,\omega_{s_1}\cdots \omega_{s_n}) \mapsto \int_\gamma \omega_{s_1}\cdots \omega_{s_n}.$$
	\end{enumerate}
	
	\begin{remark}
	The construction of the ind-motive $M$ relies on a theorem of Beilinson  \cite{goncharovmultiple, delignegoncharov} which identifies the truncation $\QQ[\pi_1(X^{\mathrm{an}},a,b)]/I^{n+1}$ with an explicit relative homology group. The isomorphism \eqref{eq: comp for pi_1} is originally due to Chen \cite{cheniterated} in the context of rational homotopy theory, and was used by Hain \cite{hainfundamentalgroup} to define the mixed Hodge structure on $M_\B$  without the use of cohomology or motives.
	\end{remark}
	
	A very useful variant of this construction, designed by Deligne  \cite{deligneP1}, replaces the basepoint $a$ and/or $b$ by a \emph{tangential basepoint at infinity}. Such an object is the datum of a point $s\in S$ together with a non-zero tangent vector $v\in (T_s\mathbb{P}^1_F)^\times$. The de Rham theory does not change (in fact, it is clear by the description above that $M_\dR$ does not depend on basepoints) but the corresponding Betti theory is concerned with paths $\gamma$ in $\mathbb{P}^1(\CC)$ that may start or end ``at infinity'', i.e., at the point $s\in S$, with prescribed outgoing or incoming velocity $v$. The extra datum of $v$ is needed to give a meaning to the possibly divergent iterated integral of a word of $1$-forms along $\gamma$ (this is called \emph{logarithmic regularization}, see the next example). 
	
	\begin{example}
	Consider $X=\PP^1_\QQ\setminus \{0,\infty\}$ and fix a tangential basepoint $(0,v)$ at $0$, with $v\in (T_0\mathbb{P}^1_\QQ)^\times= \QQ^\times$, and a classical basepoint $z\in X(\QQ)=\QQ^\times$. Even though the form $\d x/x$ is not integrable near $0$, one can make sense of its integral along a smooth path $\gamma:[0,1]\to \mathbb{P}^1(\CC)$ satisfying $\gamma(0)=0$, $\gamma'(0)=v$, $\gamma(1)=z$, and $\gamma(t)\notin\{0,\infty\}$ for $t\in (0,1]$. For this, the recipe is to consider the integral
	$$\int_{\varepsilon}^1 \frac{\d\gamma(t)}{\gamma(t)} \qquad (\varepsilon>0)$$
	and then let $\varepsilon\to 0$ while formally discarding the terms $\log(\varepsilon)$. Since $\gamma(\varepsilon)\sim \varepsilon v$, one sees that the result is a determination of $\log(z/v)$. We refer the reader to \cite{DPP} for an interpretation of this recipe in cohomological terms in the language of logarithmic geometry.
	\end{example}

	Things start to get interesting for $X=\PP^1_\QQ\setminus \{0,1,\infty\}$. For instance, $\Li_n(z)$ has a natural representation as an iterated integral on $X$, for the basepoints $0$ (with any choice of non-zero tangent vector) and $z$:
	$$\Li_n(z) = \int_{0\leq t_1\leq \cdots \leq t_n\leq 1} \frac{z\,\d t_1}{1-zt_1}\frac{\d t_2}{t_2} \cdots \frac{\d t_n}{t_n} = - \int_0^z \omega_{1}\omega_0\cdots \omega_0.$$
	This is related to our first integral representation \eqref{eq: polylog integral} by the change of variables
	$$(t_1,t_2,\ldots,t_n) = (x_1x_2\cdots x_n,x_2\cdots x_n,\ldots, x_n).$$
	Other iterated integrals for $\mathbb{P}^1_F\setminus S$ include the \emph{multiple polylogarithms} \cite{goncharovICM}
	\begin{equation}\label{eq: multiple polylog}
	\Li_{n_1,\ldots,n_r}(z_1,\ldots,z_r) = \sum_{1\leq k_1<\cdots <k_r} \frac{z_1^{k_1}\cdots z_r^{k_r}}{k_1^{n_1}\cdots k_r^{n_r}} \qquad (n_i\in\mathbb{N}_{\geq 1}, |z_i|<1).
	\end{equation}

	\subsection{Multiple zeta values and the motivic fundamental groupoid of $\mathbb{P}^1\setminus\{0,1,\infty\}$}
	
	For integers $n_1,\ldots,n_{r-1}\geq 1$ and $n_r\geq 2$, one can set $z_1=\cdots =z_r=1$ in \eqref{eq: multiple polylog}, which produces the \emph{multiple zeta values}
	$$\zeta(n_1,\ldots,n_r) = \sum_{1\leq k_1<\cdots <k_r} \frac{1}{k_1^{n_1}\cdots k_r^{n_r}}.$$
	They are iterated integrals on $\mathbb{P}^1_\QQ\setminus \{0,1,\infty\}$:
	\begin{equation}\label{eq: MZV as iterated integral}
	\zeta(n_1,\ldots,n_r) = (-1)^r \int_0^1 \underbrace{\omega_1\omega_0\omega_0\cdots \omega_0}_{n_1}\cdots\underbrace{\omega_1\omega_0\omega_0\cdots \omega_0}_{n_r},
	\end{equation}	
	where the path from $0$ to $1$ is simply the interval $[0,1]$. Therefore, they are periods of the motivic fundamental group
	\begin{equation}\label{eq: motivic fundamental group MZV}
	\mathcal{O}(\pi_1^{\operatorname{mot}}(\mathbb{P}^1_\QQ\setminus\{0,1,\infty\}, (0,v),(1,w))).
	\end{equation}
	Here the vector $v$ corresponds to $1\in \QQ^\times = (T_0\mathbb{P}^1_\QQ)^\times$, and $w$ corresponds to $-1\in\QQ^\times=(T_1\mathbb{P}^1_\QQ)^\times$. An important remark is that the tangential basepoints $(0,v)$ and $(1,w)$ are \emph{defined over $\ZZ$} in the sense that $v$ and $w$ are $\ZZ$-points of the punctured tangent spaces of $\mathbb{P}^1_\ZZ$ at $0$ and $1$ respectively, which are copies of $\mathbb{G}_{m,\ZZ}$. This is the reason why \eqref{eq: motivic fundamental group MZV} lives in the full tannakian subcategory
	$$\MT(\ZZ)\subset \MT(\QQ)$$
	of mixed Tate motives over $\ZZ$, also known as \emph{unramified} mixed Tate motives over $\QQ$ \cite{delignegoncharov}. 
	
	Note that there are only $6$ basepoints of $\mathbb{P}^1_\QQ\setminus\{0,1,\infty\}$ which are defined over $\ZZ$, namely $(0,v)$, $(0,-v)$, $(1,w)$, $(1,-w)$, and the two tangential basepoints at $\infty$ corresponding to $1$ and $-1$ in $\QQ^\times=(T_\infty\mathbb{P}^1_\QQ)^\times$. They form a torsor under the action of the symmetric group $\mathfrak{S}_3$ on $\mathbb{P}^1$, generated by the involutions $x\mapsto 1-x$ and $x\mapsto 1/x$, that globally stabilizes $\{0,1,\infty\}$ (symmetries of the cross-ratio). It is therefore natural to consider the motivic fundamental \emph{groupoid} of $\mathbb{P}^1_\QQ\setminus \{0,1,\infty\}$ relative to those $6$ basepoints. More generally, one can replace $\mathbb{P}^1_\QQ\setminus \{0,1,\infty\}$ with the moduli space
	$$\mathcal{M}_{0,n}$$
	of (smooth projective) genus zero curves with $n$ marked points, and consider its motivic fundamental groupoid relative to the finite set of tangential basepoints at infinity defined over $\ZZ$. This object has a rich (operadic) structure and plays a central role in Grothendieck--Teichm\"{u}ller theory \cite{barnatan}, a fascinating subject in itself which goes beyond the scope of this survey.

	\subsection{Volumes in hyperbolic geometry}
	
		We follow Goncharov's seminal article \cite{goncharovvolumes} (see also \cite{browndedekind} and \cite{rudenkodepth} for more recent developments). Consider an odd-dimensional hyperbolic simplex $\Sigma\subset \mathbb{H}^{2n-1}$, with $n\geq 1$. In the Klein model, 
		$$\mathbb{H}^{2n-1} = \{(x_1,\ldots,x_{2n-1})\in\RR^{2n-1} \; | \; x_1^2+\cdots +x_{2n-1}^2<1\}$$
		is the unit ball, and $\Sigma$ is bounded by a union $D=D_0\cup \cdots \cup D_{2n-1}$ of hyperplanes in $\RR^{2n-1}$. Up to complexifying and projectifying we can view everything in complex projective space $\PP^{2n-1}(\CC)$, where the boundary of $\mathbb{H}^{2n-1}$ is given by the real points of the quadric
		$$Q=\{q=0\} \quad \mbox{ with } \quad q = x_0^2 - x_1^2 - \cdots - x_{2n-2}^2 - x_{2n-1}^2.$$
		The hyperbolic volume of $\Sigma$ is given by the formula
		\begin{equation}\label{eq: hyperbolic volume}
		\operatorname{vol}(\Sigma) =  \int_\Sigma\omega_Q,
		\end{equation}
		where
		$$\omega_Q= \frac{\d x_1\cdots \d x_{2n-1}}{\left(1-x_1^2-\cdots -x_{2n-1}^2\right)^n}$$
		extends to an algebraic $(2n-1)$-form on $\PP^{2n-1}\setminus Q$. It is therefore clear that \eqref{eq: hyperbolic volume} is a period of the relative cohomology group
			\begin{equation}\label{eq: relative cohomology hyperbolic volume}
			\H^{2n-1}(\PP^{2n-1} \setminus Q , D \setminus Q\cap D).
			\end{equation}
			One can then consider such cohomology groups for any smooth quadric $Q$ and any union of hyperplanes $D$. If $Q$ and $D$ are defined over a number field $F$ and if $Q$ splits over $F$, then \eqref{eq: relative cohomology hyperbolic volume} lifts to a mixed Tate motive over $F$. This is used in \cite{browndedekind} to give formulas for special values of the Dedekind zeta function of a totally real number field as in \eqref{eq: general formula for zeta F n}.
			
	\subsection{Feynman integrals and motives}
	
		In modern particle physics, Feynman integrals are an essential tool to predict the interactions between particules at high energy (for instance in a particle collider such as the Large Hadron Collider). To a \emph{Feynman graph} $\Gamma$ (a graph in the usual mathematical sense equipped with certain physical decorations) with $n$ edges and first Betti number $h$ one associates the $(n-1)$-fold integral
	\begin{equation}\label{eq: feynman integral}
	I_\Gamma = \int_{\sigma} \frac{\Phi_\Gamma^{n-2(h+1)}}{\Xi_\Gamma^{n-2h}}\sum_{i=1}^{n}(-1)^ix_i\, \d x_1\wedge\cdots\wedge\widehat{\d x_i}\wedge \cdots \wedge \d x_n,
	\end{equation}
	where $\Phi_\Gamma$, $\Xi_\Gamma$ are the \emph{graph polynomials} (homogeneous of degree $h$ and $h+1$ respectively), and $\sigma$ is the locus in $\mathbb{P}^{n-1}(\RR)$ where $x_1,\ldots,x_n\geq 0$. The cohomological and motivic study of Feynman integrals was launched by Bloch--Esnault--Kreimer \cite{blochesnaultkreimer} and later developed by Brown \cite{browncosmic}. One finds interesting mixed Tate motives in this way, although the motives underlying \eqref{eq: feynman integral} are rarely mixed Tate.
		
\section{The tannakian formalism for mixed Tate motives: tools}\label{sec: tannakian}

	We now survey the tools provided by the tannakian formalism for mixed Tate motives, which were first advertised and used by Beilinson--Deligne \cite{beilinsondeligne} and Goncharov \cite{goncharovseattle}, and further developed by Brown \cite{brownnotes}. Examples will follow in the next section.

	\subsection{Motivic Galois groups}
	
	Let $F$ be a field for which the Beilinson--Soul\'{e} vanishing conjecture holds. The tensor product of $\DM(F)$ induces the structure of a tannakian category on $\MT(F)$, for which \eqref{eq: fiber functor MT} is a fiber functor. We let
	$$G_{\MT(F)} := \underline{\operatorname{Aut}}^\otimes(\omega)$$
	denote the corresponding Tannaka group. It is an affine group scheme over $\QQ$ and $\omega$ induces an equivalence of categories
	\begin{equation}\label{eq: MT equivalent to Rep}
	\MT(F)\stackrel{\sim}{\longrightarrow} \mathsf{Rep}(G_{\MT(F)}).
	\end{equation}
	
	\begin{definition}
	We call $G_{\MT(F)}$ the \emph{motivic Galois group} of $\MT(F)$.
	\end{definition}
	
	The action of $G_{\MT(F)}$ on $\omega(\QQ(-1))=\QQ$ induces a morphism of group schemes $G_{\MT(F)}\to \mathbb{G}_m$ whose kernel is denoted by $U_{\MT(F)}$, so that we have a short exact sequence of group schemes
	$$0\longrightarrow U_{\MT(F)} \longrightarrow G_{\MT(F)} \longrightarrow \mathbb{G}_m\longrightarrow 0. $$
	A first important remark is that this short exact sequence splits because the fiber functor $\omega$ is naturally $\mathbb{Z}$-graded. More precisely, the natural splitting
	 \begin{equation}\label{eq: tau splitting}
	 \tau\colon \mathbb{G}_m\to G_{\MT(F)}
	 \end{equation} 
	 sends an invertible scalar $\lambda$ to the automorphism of $\omega$ which acts by $\lambda^n$ on the $n$-th component of the direct sum in \eqref{eq: fiber functor MT}. This implies that we have a semidirect product decomposition:
	\begin{equation}\label{eq: semidirect product motivic Galois}
	G_{\MT(F)} \simeq \mathbb{G}_m \ltimes U_{\MT(F)}.
	\end{equation}
	A second important remark is that $U_{\MT(F)}$ is a pro-unipotent group scheme, i.e., a limit of unipotent algebraic groups. Indeed, for an object $M$ of $\MT(F)$, the action of $U_{\MT(F)}$ on $\omega(M)$ respects the filtration by the subspaces $\omega(\W_{2n}M)$ and acts trivially on the successive quotients 
	$\omega(\gr_{2n}^\W M)$ because $\gr_{2n}^\W M$ is a direct sum of tensor powers of $\QQ(-1)$. Therefore,  after choosing a basis of $\omega(M)$ compatible with the grading, the action of $U_{\MT(F)}$ is given by unipotent matrices.
	
	\subsection{The motivic Hopf algebra and Lie coalgebra}
	
	As is usual with Tannaka groups, exhibiting \emph{points} for them is difficult, whereas it is easy to produce and manipulate \emph{functions} on them. This is already the case for the more classical Galois groups: only two elements of $\mathrm{Gal}(\overline{\QQ}/\QQ)$ are known (the identity and complex conjugation), whereas every quadratic extension $K$ of $\QQ$ gives rise to a character $\mathrm{Gal}(\overline{\QQ}/\QQ)\to \mathrm{Gal}(K/\QQ)\simeq \{\pm 1\}$.
	
	This suggests that the more concrete object to consider is the Hopf algebra of functions on the motivic Galois group $G_{\MT(F)}$. Since the latter is completely determined by $U_{\MT(F)}$ with its action of $\mathbb{G}_m$ by \eqref{eq: semidirect product motivic Galois}, we make the following definition instead.
	
	\begin{definition}
	The \emph{motivic Hopf algebra} of $\MT(F)$ is the algebra of functions on the pro-unipotent group scheme $U_{\MT(F)}$:
	$$\mathcal{H}(F):= \mathcal{O}(U_{\MT(F)}).$$
	\end{definition}
	
	It is a Hopf algebra over $\QQ$. The action of $\mathbb{G}_m$ on $\mathcal{O}(U_{\MT(F)})$ translates as a grading on $\mathcal{H}(F)$, and by \eqref{eq: MT equivalent to Rep} and \eqref{eq: semidirect product motivic Galois} one has that $\MT(F)$ is equivalent to the category of graded comodules over $\mathcal{H}(F)$:
	$$\MT(F)\stackrel{\sim}{\longrightarrow} \mathsf{grComod}(\mathcal{H}(F)).$$

	If $F$ is a number field then \eqref{eq: ext in MT F number field} implies that $\mathcal{H}(F)$ is a cofree Hopf algebra cogenerated in degree $n\geq 1$ by a copy of $\K_{2n-1}(F)_\QQ$. Therefore it is non-negatively graded and connected, i.e.,
	$$\mathcal{H}(F) = \bigoplus_{n\geq 0} \mathcal{H}_n(F) \qquad \mbox{ with }\; \mathcal{H}_0(F)=\QQ,$$	
	 and one can ``compute'' the rational $\K$-theory of $F$ as the space of primitive elements of $\mathcal{H}(F)$:
	\begin{equation}\label{eq: K theory as primitive elements}
	\K_{2n-1}(F)_\QQ \simeq \ker\left(\overline{\Delta}:\mathcal{H}_n(F)\longrightarrow \bigoplus_{k=1}^{n-1}\mathcal{H}_k(F)\otimes \mathcal{H}_{n-k}(F)\right),
	\end{equation}
	where $\overline{\Delta}(x):=\Delta(x)-1\otimes x-x\otimes 1$ is the reduced coproduct.
	One can also consider the space of indecomposables
	$$\mathcal{C}(F) := \mathcal{H}_{>0}(F)/\mathcal{H}_{>0}(F)\mathcal{H}_{>0}(F),$$
	where $\mathcal{H}_{>0}(F):=\bigoplus_{n>0}\mathcal{H}_n(F)$. The coproduct of $\mathcal{H}(F)$ induces a \emph{cobracket}
	$$\delta:\mathcal{C}(F)\longrightarrow \Lambda^2\mathcal{C}(F)$$
	which gives $\mathcal{C}(F)$ the structure of a Lie coalgebra. It is sometimes called the \emph{motivic Lie coalgebra} and contains as much information as the motivic Hopf algebra, which can be recovered as its universal coenveloping coalgebra. In particular, we have the following variant of \eqref{eq: K theory as primitive elements}:
	\begin{equation}\label{eq: K theory inside Lie coalgebra}\K_{2n-1}(F)_\QQ \simeq \ker\left(\delta : \mathcal{C}_n(F)\longrightarrow (\Lambda^2\mathcal{C}(F))_n \right).
	\end{equation}
	
	\subsection{Matrix coefficients} 
	
		The motivic Hopf algebra $\mathcal{H}(F)$ is a rather concrete object. By the general tannakian formalism, the Hopf algebra
		$$\mathcal{O}(G_{\MT(F)})$$
		is spanned by matrix coefficients
		\begin{equation}\label{eq: matrix coefficient} 
		(M,\varphi,v)
		\end{equation}
		with $M$ an object of $\MT(F)$, $v\in \omega(M)$ and $\varphi\in \omega(M)^\vee$. The function on $G_{\MT(F)}$ corresponding to \eqref{eq: matrix coefficient} is 
		$$g\mapsto \langle \varphi,g\cdot v \rangle,$$
		where $\langle \cdot,\cdot\rangle$ is the duality pairing. The $\QQ$-linear relations among matrix coefficients are spanned by the obvious relations
		$$(M,\psi,f(v)) = (N,f^\vee(\psi),v)$$
		for  morphisms $f\colon M\to N$ in $\MT(F)$ and $v\in \omega(M)$, $\psi\in\omega(N)^\vee$. The coproduct of $\mathcal{O}(G_{\MT(F)})$ is computed on matrix coefficients by the formula
		$$\Delta(M,\varphi,v) = \sum_{i} (M,\varphi,e_i)\otimes (M,e_i^\vee,v),$$
		where $(e_i)$ is any basis of $\omega(M)$ and $(e_i^\vee)$ denotes the dual basis.
		
		\begin{example}
		We define the \emph{Lefschetz element} of $\mathcal{O}(G_{\MT(F)})$ as 
		$$\mathbb{L} := (\QQ(-1),1^\vee,1)$$
		where $1$ denotes the canonical basis element of $\omega(\QQ(-1))=\QQ$. It is a group-like element of the Hopf algebra $\mathcal{O}(G_{\MT(F)})$ in the sense that
		$$\Delta(\mathbb{L})=\mathbb{L}\otimes\mathbb{L}.$$
		\end{example}
		
		The motivic Hopf algebra $\mathcal{H}(F)$ is the quotient of $\mathcal{O}(G_{\MT(F)})$ by the (Hopf) ideal generated by $(\mathbb{L}-1)$, which has the effect of ``trivializing'' Tate twists. The space $\mathcal{H}_n(F)$ of homogeneous elements of degree $n$ is spanned by the matrix coefficients \eqref{eq: matrix coefficient} for which $\varphi:\gr_0^\W M\to \QQ(0)$ and $v:\QQ(-n)\to \gr_{2n}^\W M$ .
		
		\begin{example}
		Let $x\in F\setminus \{0,1\}$. By the short exact sequence \eqref{eq: short exact sequence Kummer motive}  there are canonical identifications $\gr_0^\W \mathcal{K}_x\simeq \QQ(0)$ and $\gr_2^\W\mathcal{K}_x\simeq \QQ(-1)$. We denote the corresponding matrix coefficient by 
		$$\log^{\mathcal{H}}(x) \; \in \mathcal{H}_1(F).$$
		It is a primitive element of the Hopf algebra $\mathcal{H}(F)$ in the sense that
		$$\Delta(\log^{\mathcal{H}}(x)) = 1\otimes \log^{\mathcal{H}}(x) + \log^{\mathcal{H}}(x)\otimes 1.$$
		The notation is justified by the fact that $\log^{\mathcal{H}}$ satisfies the functional equation of the logarithm:
		\begin{equation}\label{eq: log H functional equation}
		\log^{\mathcal{H}}(xy)=\log^{\mathcal{H}}(x) +\log^{\mathcal{H}}(y).
		\end{equation}
		\end{example}

	\subsection{Regulators as single-valued periods}\label{subsec: regulators sv}
	
	Let $F$ be a number field and $\sigma\colon F\to \CC$ be an embedding. The single-valued period isomorphism \eqref{eq: sv sigma general} induces, thanks to \eqref{eq: omega dR vs omega canonical}, a $\CC$-linear isomorphism
	$$\operatorname{s}_\sigma:\omega\otimes_\QQ\CC\stackrel{\sim}{\longrightarrow} \omega\otimes_\QQ\CC,$$
	i.e., a  complex point of the motivic Galois group $G_{\MT(F)}$. By the computation of Example \ref{ex: sv pi}, its image by the map $G_{\MT(F)}\to \mathbb{G}_m$ is $-1$, therefore it does not live in the pro-unipotent group scheme $U_{\MT(F)}$. It is however possible to correct it and set
	$$\operatorname{sv}_\sigma^{\mathcal{H}} := \tau(-1)\operatorname{s}_\sigma,$$
	where $\tau$ is the splitting \eqref{eq: tau splitting}, which defines a complex point of $U_{\MT(F)}$, or equivalently a morphism of algebras
	$$\operatorname{sv}_\sigma^{\mathcal{H}}:\mathcal{H}(F)\to \CC,$$
	sometimes called the \emph{single-valued period map} (for mixed Tate motives over $F$).
	
	Now since $U_{\MT(F)}$ is pro-unipotent one can take the logarithm of $\operatorname{sv}_\sigma^{\mathcal{H}}$, which lives in its Lie algebra, i.e.,  is a complex-valued function on the motivic Lie coalgebra. We then define
	$$\operatorname{sv}_\sigma^{\mathcal{C}}:= \frac{1}{2} \log(\operatorname{sv}_\sigma^{\mathcal{H}}):\mathcal{C}(F)\to \CC$$
	and call it the \emph{single-valued Lie-period map}. In view of the definition of $\mathcal{C}(F)$ as the space of indecomposables in $\mathcal{H}(F)$, the map $\operatorname{sv}_\sigma^{\mathcal{C}}$ ``kills products'', which makes it better behaved than $\operatorname{sv}_\sigma^{\mathcal{H}}$ in certain situations.

	Concretely, an element of $\mathcal{C}_n(F)$ is an equivalence class of a matrix coefficient $(M,\varphi,v)$ with $M\in\MT(F)$, $\varphi:\gr_0^\W M\to \QQ(0)$, and $v:\QQ(-n)\to \gr_{2n}^\W M$, considered modulo products. Let $P$ be an $(\omega,\omega_\B)$-period matrix for $M$, i.e. a matrix for the comparison isomorphism
	$$\omega(M)\otimes_\QQ\CC\simeq \omega_\dR(M)\otimes_{F,\sigma}\CC\stackrel{\sim}{\To} \omega_{\B,\sigma}(M)\otimes_\QQ\CC,$$
	relative to a graded basis of $\omega(M)$. Let $D$ be the basis of $\tau(-1)$ acting on $\omega(M)$, i.e., the diagonal matrix with entry $(-1)^k$ in weight $2k$. The product $D\overline{P}^{\,-1}P$ is unipotent, and the single-valued Lie-period map is then computed as
	$$\operatorname{sv}_\sigma^{\mathcal{C}}(M,\varphi,v) =  \frac{1}{2}\langle \varphi,\log(D\overline{P}^{\,-1}P)v\rangle.$$
	
	\begin{example}
	Let $M\in\MT(F)$ be an extension
	$$0\To \QQ(0)\To M\To \QQ(-n)\To 0,$$
	and let $\xi=(M,\varphi,v)\in C_n(F)$ denote the corresponding matrix coefficient. Write the period matrix of $M$ in the form
	$$P=\setlength{\arraycolsep}{3pt}\def\arraystretch{1.3}
		\left(\begin{matrix}  1 & \alpha \\ 0 & (2\pi\i)^n\end{matrix}\right).$$	
	Then we have
	$$\log(D\overline{P}^{\, -1}P) = \setlength{\arraycolsep}{4pt}\def\arraystretch{1.5}
		\left(\begin{matrix}  0 & \alpha-(-1)^n\overline{\alpha} \\ 0 & 0\end{matrix}\right) $$
	and therefore
	$$\operatorname{sv}_\sigma^{\mathcal{C}}(\xi) = \begin{cases} \operatorname{Re}(\alpha) & \mbox{ if } n \mbox{ is odd} \\ \i\operatorname{Im}(\alpha) & \mbox{ if } n \mbox{ is even.}\end{cases}$$
	For instance, in the case $n=1$, we have $\alpha=\log(\sigma(x))$ and therefore
	\begin{equation}\label{eq: sv of log}
	\operatorname{sv}^{\mathcal{C}}_\sigma(\log^{\mathcal{H}}(x)) = \operatorname{Re}(\log(\sigma(x))) = \log|\sigma(x)|.
	\end{equation}
	\end{example}
	
	This example gives a concrete way of computing the Hodge regulator via the following commutative diagram, where the horizontal inclusion is induced by \eqref{eq: K theory inside Lie coalgebra}:
	$$\xymatrixrowsep{1.2cm}\diagram{
	\K_{2n-1}(F)\; \ar@{^(->}[rr] \ar[rd]_{\varpi_n^{(\sigma)}} && \;C_n(F) \ar[dl]^{\operatorname{sv}_\sigma^{\mathcal{C}}} \\
	 & \CC/(2\pi\i)^n\QQ &
	}$$
	
\section{The tannakian formalism for mixed Tate motives: examples}\label{sec: tannakian examples}

	Let $F$ be a number field. We use our examples from \S\ref{sec: families MT} to produce elements of the motivic Hopf algebra (resp. the motivic Lie coalgebra) of $\MT(F)$ and compute their coproduct (resp. their cobracket).

	\subsection{Motivic polylogarithms}
	
		For $x\in F\setminus \{0,1\}$ and $n\geq 1$, recall the polylogarithm motive $\mathcal{L}_n^{(x)}$, for which we have canonical isomorphisms
		$$\gr_0^\W \mathcal{L}_x^{(n)} \simeq \QQ(0) \quad \mbox{ and } \quad \gr_{2n}^\W \mathcal{L}_x^{(n)}\simeq \QQ(-n),$$
		induced by \eqref{eq: polylog as extension} and \eqref{eq: polylog inductive}. We denote by 
		$$\Li_n^{\mathcal{H}}(x) \in\mathcal{H}_n(F)$$
		the corresponding matrix coefficient, and by
		$$\Li_n^{\mathcal{C}}(x)\in \mathcal{C}_n(F)$$
		its image modulo products. They are sometimes referred to as \emph{motivic polylogarithms}. As explained in \S\ref{subsec: sv polylog}, we have
		$$\operatorname{sv}_\sigma^{\mathcal{C}}(\Li_n^{\mathcal{C}}(x)) = \operatorname{P}_n(\sigma(x)).$$
		
		From the structure \eqref{eq: polylog as extension} and \eqref{eq: polylog inductive} of the polylogarithm motive one easily derives the coproduct formula:
		$$\Delta(\Li_n^{\mathcal{H}}(x)) = \sum_{k=0}^{n-1} \Li_{n-k}^{\mathcal{H}}(x)\otimes \frac{(\log^{\mathcal{H}}(x))^k}{k!} + 1\otimes \Li_n^{\mathcal{H}}(x).$$
		Modulo products, we get an even simpler cobracket formula:
		\begin{equation}\label{eq: cobracket polylog}
		\delta(\Li_n^{\mathcal{C}}(x)) = \Li_{n-1}^{\mathcal{C}}(x)\wedge \log^{\mathcal{C}}(x).
		\end{equation}
		This formula is reminiscent of the differential equation \eqref{eq: diff eq Li n} for the classical polylogarithms, written as $\d\Li_n(z) = \Li_{n-1}(z)\operatorname{dlog}(z)$. This is not a coincidence and can be explained by the general compatibility between the motivic coproduct and the Gauss--Manin connection for motivic periods over a base \cite[\S 7]{brownnotes}.
		
		\emph{Zagier's polylogarithm conjecture} predicts that motivic polylogarithms are enough to capture all primitive elements in the motivic Lie coalgebra, i.e., to span the algebraic $\K$-theory of number fields by \eqref{eq: K theory inside Lie coalgebra}.
		
		\begin{conjecture}\label{conj: zagier}
		The space of primitive elements in $\mathcal{C}_n(F)$ is spanned by $\QQ$-linear combinations of elements $\Li_n^{\mathcal{C}}(x)$, for $x\in F\setminus \{0,1\}$.
		\end{conjecture}
		
		Under the equality of regulators \eqref{eq: equality of regulators}, it implies Zagier's original conjecture \cite{zagierconjecture}, which is an expression \eqref{eq: general formula for zeta F n} where each $\alpha_j^{(\sigma)}$ is of the form
		$$\sum_k a_k\operatorname{P}_n(\sigma(x_{k}))$$
		for some $a_k\in\QQ$ and $x_k\in F\setminus\{0,1\}$. This has been proved by Zagier  for $n=2$, Goncharov \cite{goncharovgeometry} for $n=3$, and Goncharov and Rudenko \cite{goncharovrudenko} for $n=4$. We refer the reader to the survey article \cite{dupontbourbaki} for more details on Zagier's conjecture and its motivic aspects.
		
		\begin{remark}
		The cobracket formula \eqref{eq: cobracket polylog} is so simple that it sometimes allows one to produce many primitive elements in $\mathcal{C}_n(F)$.  For instance, if $x$ is a root of unity then $\log^{\mathcal{C}}(x)=0$ and therefore $\operatorname{Li}_n^{\mathcal{C}}(x)$ is primitive. One can prove that those elements $\operatorname{Li}_n^{\mathcal{C}}(x)$ span the space of primitive elements if $F$ is a cyclotomic field, which settles Zagier's conjecture in this case. For another example borrowed from \cite{zagiergangl}, take $F=\QQ(\sqrt{-7})$, $n=2$, and consider the element
		$$\xi = 2 \operatorname{Li}_2^{\mathcal{C}}(y) + \operatorname{Li}_2^{\mathcal{C}}(y') \quad \mbox{ for } y=\frac{1+\sqrt{-7}}{2} \mbox{ and } y'=\frac{-1+\sqrt{-7}}{4}.  $$
		Since $\operatorname{Li}_1^{\mathcal{C}}(x)= -\log^{\mathcal{C}}(1-x)$, one easily proves using \eqref{eq: cobracket polylog} and \eqref{eq: log H functional equation} that $\delta(\xi)=0$. (Hint: use $1-y=-2y'$ and $1-y'=-y^2$.)
		\end{remark}
	
	\subsection{Motivic iterated integrals}
	
	Using motivic fundamental group(oid)s one can produce ``motivic'' versions of iterated integrals on the punctured projective line. Following \cite{goncharovgaloissymmetries}, we write
	$$\mathbb{I}(a_0;a_1,\ldots,a_n;a_{n+1}) := \int_{a_0}^{a_{n+1}} \omega_{a_1}\cdots \omega_{a_n},$$
	for elements $a_0,\ldots,a_{n+1}\in F$ (one potentially needs to specify tangential basepoints at $a_0$ and $a_{n+1}$), and define
	\begin{equation}\label{eq: motivic iterated integral}
	\mathbb{I}^{\mathcal{H}}(a_0;a_1,\ldots,a_n;a_{n+1}) \in \mathcal{H}_n(F)
	\end{equation}
	to be the matrix coefficient corresponding to the mixed Tate motive
	\begin{equation}\label{eq: some motivic fundamental group}
	M = \mathcal{O}(\pi_1^{\operatorname{mot}}(\mathbb{P}^1_F\setminus\{a_1,\ldots,a_n\},a_0,a_{n+1}))
	\end{equation}
	with $\varphi:\gr_0^\W M\stackrel{\sim}{\to}\QQ(0)$ the obvious identification, and $v:\QQ(-n)\to \gr_{2n}^W M$ given by the word $\omega_{a_1}\cdots \omega_{a_n}$. The coproduct of the elements \eqref{eq: motivic iterated integral} is constrained by the structure of motivic fundamental group(oid)s and given by \emph{Goncharov's formula}
	\begin{equation}\label{eq: goncharov formula}
	\begin{split}
	\Delta(\mathbb{I}^{\mathcal{H}}(a_0;a_1,\ldots,a_n;a_{n+1}))  & = \sum_{\substack{0\leq k\leq n\\0=i_0<i_1<\cdots<i_k<i_{k+1}=n+1}}  \\
	\left(\prod_{s=0}^k \mathbb{I}^{\mathcal{H}}(a_{i_s};a_{i_s+1},\ldots, a_{i_{s+1}-1};a_{i_{s+1}})  \right)\;&\otimes \; \mathbb{I}^{\mathcal{H}}(a_0;a_{i_1},\ldots, a_{i_k};a_{n+1}).
	\end{split}
	\end{equation}
 	Its cobracket version (i.e., modulo products) simplifies as
	\begin{equation*}
			\begin{split}
			\delta &(\mathbb{I}^{\mathcal{C}}(a_0;a_1,\ldots,a_n;a_{n+1}))  \\
			& = \sum_{0\leq i<j\leq n} \mathbb{I}^{\mathcal{C}}(a_0;a_1,\ldots,a_i,a_{j+1},\ldots,a_n;a_{n+1})
			\wedge \mathbb{I}^{\mathcal{C}}(a_i;a_{i+1},\ldots,a_{j-1};a_j)\ .
			\end{split}
			\end{equation*}
			
		Goncharov \cite{goncharovICM} conjectures that iterated integrals are enough to understand the structure of mixed Tate motives.
			
		\begin{conjecture}\label{conj: universality}
		The motivic iterated integrals \eqref{eq: motivic iterated integral}, for all choices of $a_i\in F$, span the motivic Hopf algebra $\mathcal{H}(F)$.
		\end{conjecture}
		
		\begin{remark}
		By the tannakian dictionary, this would translate as the fact that motivic fundamental groups \eqref{eq: some motivic fundamental group} generate $\MT(F)$ as a tannakian category.
		\end{remark}
		
		Prompted by Conjecture \ref{conj: universality}, one may dream of giving a complete ``combinatorial'' description of the motivic Hopf algebra $\mathcal{H}(F)$ by generators and relations. This idea has been pursued by Goncharov who gave such conjectural descriptions in low weight and proved that they are compatible with known descriptions of the algebraic $\K$-theory of fields. We refer the reader to \cite{dupontbourbaki} for an introduction to this circle of ideas, which we call the \emph{Goncharov program}.
		
		\begin{remark}
		Zagier's conjecture (Conjecture \ref{conj: zagier}) does not imply Conjecture \ref{conj: universality} because it only gives information on the primitive elements in $\mathcal{H}(F)$. Conversely, Conjecture \ref{conj: universality} does not imply Zagier's conjecture because the latter predicts that a very specific family of motivic iterated integrals, namely the motivic versions of classical polylogarithms, span the primitives of $\mathcal{H}(F)$. 
		\end{remark}

	\subsection{Motivic multiple zeta values}
	
	As a special case of the construction of motivic iterated integrals, using \eqref{eq: MZV as iterated integral}, one gets motivic multiple zeta values 
	$$\zeta^{\mathcal{H}}(n_1,\ldots,n_r) \in \mathcal{H}_n(\ZZ) \quad \mbox{ with } n=n_1+\cdots +n_r,$$
	where $\mathcal{H}(\ZZ)$ denotes the motivic Hopf algebra of the category $\MT(\ZZ)\subset \MT(\QQ)$. Brown \cite{brownMTZ, brownICM}	 defined certain variants  that we denote by
	\begin{equation}\label{eq: true motivic MZVs}
	\zeta^{\mathcal{P}}(n_1,\ldots,n_r)\in\mathcal{P}_n(\ZZ).
	\end{equation}
	They live in the \emph{algebra of motivic periods} for $\MT(\ZZ)$, defined as 
	$$\mathcal{P}(\ZZ) := \mathcal{O}(\underline{\operatorname{Isom}}^\otimes(\omega_\dR,\omega_\B)),$$
	where $\omega_\dR,\omega_\B:\MT(\ZZ)\to \mathsf{Vect}_\QQ$ are the de Rham and Betti fiber functors. The comparison isomorphism \eqref{eq: comparison isomorphism} induces a morphism of algebras
	\begin{equation}\label{eq: period map MT Z}
	\operatorname{per}:\mathcal{P}(\ZZ) \To \CC,
	\end{equation}
	called the \emph{period map}. It allows one to recover multiple zeta values from Brown's motivic versions \eqref{eq: true motivic MZVs}:
	$$\operatorname{per}(\zeta^{\mathcal{P}}(n_1,\ldots,n_r)) = \zeta(n_1,\ldots,n_r).$$
	The Hopf algebra structure on $\mathcal{H}(\ZZ)$ is replaced by a \emph{coaction} \cite{brownnotes}
	$$\mathcal{P}(\ZZ) \To \mathcal{P}(\ZZ)\otimes\mathcal{H}(\ZZ),$$
	which can be computed on \eqref{eq: true motivic MZVs} by a variant of Goncharov's formula \eqref{eq: goncharov formula}. These tools are the main inputs of the proof of the following important theorem of Brown \cite{brownMTZ} (previously known as the Deligne--Ihara conjecture).
	
	\begin{theorem}
	The object $\mathcal{O}(\pi_1^{\operatorname{mot}}(\mathbb{P}^1\setminus\{0,1,\infty\},(0,v),(1,w)))$ generates $\MT(\ZZ)$ as a tannakian category.
	\end{theorem}
	
	The proof of this theorem actually gives a more detailed description of the structure of mixed Tate motives over $\ZZ$. Using the period map \eqref{eq: period map MT Z}, one gets the following theorem (conjectured by Hoffman \cite{hoffman}), whose statement does not involve mixed Tate motives at all, but whose only proof (at the time of writing) requires the help of the motivic machinery.

	\begin{theorem}
	Every multiple zeta value can be expressed as a $\QQ$-linear combination of multiple zeta values $\zeta(n_1,\ldots,n_r)$ for which all $n_i\in\{2,3\}$.
	\end{theorem}
	
	We refer the reader to the survey article \cite{dupontxups} and the monograph \cite{burgosfresan} for more details on multiple zeta values, their motivic versions, and Brown's theorem.

	\subsection{Motivic volumes in hyperbolic geometry}
	
		Using the mixed Tate motives which lift the relative cohomology groups \eqref{eq: relative cohomology hyperbolic volume}, Goncharov \cite{goncharovvolumes} defines a motivic version of the volume of a hyperbolic simplex $\Sigma$ in $\mathbb{H}^{2n-1}$ whose faces (in Klein's model) are defined over a number field $F$,
		\begin{equation}\label{eq: motivic volume Sigma}
		\operatorname{vol}^{\mathcal{H}}(\Sigma) \in \mathcal{H}_n(F).
		\end{equation}
		Its coproduct
		$$\Delta(\operatorname{vol}^{\mathcal{H}}(\Sigma)) \in\mathcal{H}(F)\otimes\mathcal{H}(F)$$
		is related to the Dehn invariant of $\Sigma$, which was used to solve Hilbert's third problem by proving that two polyhedra with equal volume cannot always be transformed into one another by cutting and pasting.
		
		By triangulating a given hyperbolic $(2n-1)$-manifold $M$ of finite volume, Goncharov defines a motivic version
		$$\operatorname{vol}^{\mathcal{H}}(M)\in\mathcal{H}_n(F)$$
		of its volume as a sum of elements \eqref{eq: motivic volume Sigma}, and proves that it gives rise to a primitive element in the motivic Hopf algebra. This builds a bridge between volumes of hyperbolic manifolds on the one hand, and regulators and special values of Dedekind zeta functions on the other hand.

\vspace{.3cm}

\bibliographystyle{alpha}
\bibliography{biblio}

\end{document}